\documentclass[a4paper,12pt]{amsart}
\usepackage{mathrsfs}
\usepackage{amsfonts}
\usepackage{amssymb}
\usepackage{ifthen}
\usepackage{graphicx}
\usepackage[usenames]{color}

\nonstopmode \numberwithin{equation}{section}
\newtheorem{thm}{Theorem}[section]
\newtheorem{cor}{Corollary}[section]
\newtheorem{lem}{Lemma}[section]
\newtheorem{proposition}{Proposition}[section]
\newtheorem{claim}{Claim}
\newtheorem{conj}{Conjecture}[section]

\theoremstyle{definition}
\newtheorem{defn}{Definition}[section]
\newtheorem{case}{Case}
\newtheorem{subcase}{Subcase}

\newtheorem{examp}{Example}
\newtheorem{prob}{Problem}
\newtheorem{ques}{Question}
\newtheorem{rem}{Remark}[section]

\newcounter {own}
\def\theown {\thesection       .\arabic{own}}

\newenvironment{pf}[1][]{%
 \vskip 3mm
 \noindent
 \ifthenelse{\equal{#1}{}}%
  {{\slshape Proof. }}%
  {{\slshape #1.} }%
 }%
{\qed\bigskip}

\newcounter{alphabet}

\newenvironment{Thm}[1][]{\refstepcounter{alphabet}%
\bigskip%
\noindent%
{\bf Theorem \Alph{alphabet}}%
\ifthenelse{\equal{#1}{}}{}{ (#1)}%
{\bf .} \itshape}{\vskip 8pt}

\newcommand{\ID}{{\mathbb D}}

\def\be{\begin{equation}}
\def\ee{\end{equation}}

\newcommand{\bee}{\begin{enumerate}}
\newcommand{\eee}{\end{enumerate}}

\newcommand{\blem}{\begin{lem}}
\newcommand{\elem}{\end{lem}}
\newcommand{\bthm}{\begin{thm}}
\newcommand{\ethm}{\end{thm}}
\newcommand{\bcor}{\begin{cor}}
\newcommand{\ecor}{\end{cor}}
\newcommand{\beg}{\begin{examp}}
\newcommand{\eeg}{\end{examp}}
\newcommand{\begs}{\begin{examples}}
\newcommand{\eegs}{\end{examples}}
\newcommand{\bdefe}{\begin{defn}}
\newcommand{\edefe}{\end{defn}}
\newcommand{\bprob}{\begin{prob}}
\newcommand{\eprob}{\end{prob}}
\newcommand{\bprop}{\begin{proposition}}
\newcommand{\eprop}{\end{proposition}}
\newcommand{\bd}{\begin{definition}}
\newcommand{\ed}{\end{definition}}
\newcommand{\bques}{\begin{ques}}
\newcommand{\eques}{\end{ques}}
\newcommand{\bei}{\begin{itemize}}
\newcommand{\eei}{\end{itemize}}
\newcommand{\bca}{\begin{case}}
\newcommand{\eca}{\end{case}}
\newcommand{\bsca}{\begin{subcase}}
\newcommand{\esca}{\end{subcase}}

\newcommand{\bcl}{\begin{claim}}
\newcommand{\ecl}{\end{claim}}

\newcommand{\bcon}{\begin{conj}}
\newcommand{\econ}{\end{conj}}
\newcommand{\bcons}{\begin{conjs}}
\newcommand{\econs}{\end{conjs}}
\newcommand{\br}{\begin{rem}}
\newcommand{\er}{\end{rem}}
\newcommand{\brs}{\begin{rems}}
\newcommand{\ers}{\end{rems}}
\newcommand{\bo}{\begin{obser}}
\newcommand{\eo}{\end{obser}}
\newcommand{\bos}{\begin{obsers}}
\newcommand{\eos}{\end{obsers}}
\newcommand{\bpf}{\begin{pf}}
\newcommand{\epf}{\end{pf}}
\newcommand{\ba}{\begin{array}}
\newcommand{\ea}{\end{array}}
\newcommand{\beq}{\begin{eqnarray}}
\newcommand{\beqq}{\begin{eqnarray*}}
\newcommand{\eeq}{\end{eqnarray}}
\newcommand{\eeqq}{\end{eqnarray*}}

\begin{document}
\bibliographystyle{amsplain}
\title {Properties for ($\alpha,\beta$)-harmonic functions}

\author[J. Qiao]{Jinjing Qiao${}^\dagger$}
\address{J. Qiao, Department of Mathematics, College of Mathematics and Information Science, Hebei University,
Baoding 071002, Hebei, People's Republic of China.}
\address{Hebei Key Laboratory of Machine Learning and Computational Intelligence, College of Mathematics and Information Science, Hebei University, Baoding 071002, Hebei, People's Republic of China.}
\email{mathqiao@126.com}

\author [J. Chang]{Jiale Chang}
\address{J. Chang, Department of Mathematics, College of Mathematics and Information Science, Hebei University, Baoding 071002, Hebei, People's Republic of China.}
\email{1353955094@qq.com}

\author[A. Rasila]{Antti Rasila}
\address{A. Rasila, Department of Mathematics with computer Science, Guangdong Technion-Israel Institute of Technology, 241 Daxue Road, Shantou 515063, Guangdong, People's Republic of China.}
\address{Department of Mathematics, Technion-Israel Institute of Technology, Haifa 3200003, Israel.}
\email{antti.rasila@iki.fi; antti.rasila@gtiit.edu.cn}



\begin{abstract}
We investigate properties of ($\alpha,\beta$)-harmonic functions. First, we discuss the coefficient estimates for ($\alpha,\beta$)-harmonic functions. In particular, we obtain Heinz's inequality for ($\alpha,\beta$)-harmonic functions, propose a coefficient bound for normalized univalent ($\alpha,\beta$)-harmonic functions and prove  that this holds for the subclass that consists of starlike functions. Furthermore, by utilizing the relationship between ($\alpha,\beta$)-harmonic functions and harmonic functions, we obtain Rad\'{o}'s theorem, Koebe type covering theorems and an area theorem. Finally, we show growth estimates and distortion estimates for ($\alpha,\beta$)-harmonic functions by using the $L^p$ norms of the boundary functions.
\end{abstract}

\keywords {($\alpha,\beta$)-harmonic function, coefficient, Heinz's inequality, Rad\'{o}'s theorem, covering theorem, area theorem, growth estimate, distortion estimate.\\
$
{}^\dagger$ Corresponding author}

\maketitle
\pagestyle{myheadings} \markboth{ J. Qiao, J. Chang and A. Rasila}{($\alpha,\beta$)-harmonic functions}

\section{Introduction }

For $a \in \mathbb{C}$ and $r > 0$, we let $\mathbb{D}(a,r) :=\{z :
|z - a| < r\}$, $\mathbb{D}_{r} := \mathbb{D}(0, r)$ and $\mathbb{D} := \mathbb{D}_1$, the open unit disk. Let $\mathbb{T} =\partial\mathbb{D}$ be the boundary of $\mathbb{D}$. 
Denote by $C^{m}(\Omega)$ the set of complex-valued $m$-times continuously differentiable functions from $\Omega$ into $\mathbb{C}$, where $\Omega$ stands for an open subset of $\mathbb{C}$ and $m$ is an nonnegative integer.  In particular, $C(\Omega):=C^{0}(\Omega)$ the set  all continuous functions on $\Omega$.
For $\alpha,\beta \in \mathbb{C}$, let
$$
L_{\alpha,\beta}:=(1-|z|^{2})\Big((1-|z|^{2})\frac{\partial^{2}}{\partial z \partial \overline{z}}+\alpha
z\frac{\partial}{\partial z}+\beta \overline{z}\frac{\partial}{\partial \overline{z}}-\alpha\beta\Big)
$$
be a second order uniformly elliptic linear partial differential operator on the unit disk $\mathbb{D}$, where
$$
\frac{\partial}{\partial z}=\frac{1}{2}\Big(\frac{\partial}{\partial x}-i\frac{\partial}{\partial y}\Big)
~~\mbox{and}~~ \frac{\partial}{\partial \overline{z}}=\frac{1}{2}\Big(\frac{\partial}{\partial x}+i\frac{\partial}{\partial y} \Big).
$$
This operator was  first introduced and studied by Geller \cite{Geller1}, and subsequently  extended by Ahern et al. \cite{Ahern} to the case of $(\alpha,\beta)$-Laplace equation of the unit ball with some additional assumptions. The
 planar case has been recently investigated  in \cite{Milos,ADEL,Markus}.

We focus on the following  homogeneous equation on $\mathbb{D}$:
\be\label{eq:491}
L_{\alpha,\beta}u=0.
\ee
A real-valued or complex-valued function $u$ in ${C}^{2}(\mathbb{D})$ is said to be $(\alpha,\beta)$-harmonic if it satisfies the equation \eqref{eq:491}. An ($\alpha,\beta$)-harmonic function $u$ on $\mathbb D$ is sense-preserving  if  $|u_{z}(z)|>|u_{\overline{z}}(z)|$ for all $z\in \mathbb D$, or sense-reversing if $|u_{z}(z)|<|u_{\overline{z}}(z)|$ throughout $\mathbb D$.

We note that if $u$ is $(\alpha,\beta)$-harmonic, then $\overline{u}$ is clearly $(\beta,\alpha)$-harmonic.  An $(0,\alpha)$-harmonic function is simply an $\alpha$-harmonic function ($\alpha>-1$) (cf. \cite{Carlsson,Olo,OLO}), while an $(\frac{\alpha}{2},\frac{\alpha}{2})$-harmonic function corresponds to a real kernel $\alpha$-harmonic function ($\alpha>-1$),
as shown, for example, in \cite{David, Li, Bo, Yong}.
Observe also that a $(0,0)$-harmonic function is  harmonic  in the usual sense. See \cite{Duren} and the references therein for the properties of harmonic mappings in the complex plane.

The function $u_{\alpha,\beta}$ 
 defined by
\be\label{eq:4203}
u_{\alpha,\beta}(z)=\frac{(1-|z|^{2})^{\alpha +\beta +1}}{(1-z)^{\alpha +1}(1-\overline{z})^{\beta
+1}}~(z\in \mathbb{D})
\ee
plays an important role in the theory of $(\alpha,\beta)$-harmonic functions.
In the following, we assume that $\alpha,\beta \in \mathbb{R}\backslash \mathbb{Z}^{-}$ satisfy $\alpha+\beta>-1$ (see \cite[Remark 6.5]{Markus} for the reason for this constraint), where $\mathbb{R}$ is the set of real numbers and $\mathbb{Z}^{-}$ is the set of negative integers. In \cite{Markus}, Klintborg and Olofsson showed that the $(\alpha,\beta)$-harmonic Poisson kernel can be defined by
$$
P_{\alpha,\beta}(z\overline{\zeta})=c_{\alpha,\beta}u_{\alpha,\beta}(z\overline{\zeta}),~~z\in \mathbb{D}, ~\zeta\in\mathbb{T},
$$
where the normalizing constant $c_{\alpha,\beta}$ is given by
$$
c_{\alpha,\beta}=\frac{\Gamma(\alpha +1)\Gamma(\beta +1)}{\Gamma(\alpha +\beta+1)}.
$$

\bdefe
For a complex-valued function $f \in L^{1}(\mathbb{T})$, the $(\alpha,\beta)$-Poisson integral of $f$ is defined by
\be\label{eq:497}
u(z)=P_{\alpha,\beta}[f](z)=\int_{\mathbb{T}}P_{\alpha,\beta}(z\overline{\zeta})f(\zeta)~dm(\zeta),~~z\in \mathbb{D}.
\ee
\edefe

It is easy to verify that for each $f \in L^{1}(\mathbb{T})$, the function $P_{\alpha,\beta}[f](z)$ is $(\alpha,\beta)$-harmonic on $\mathbb D$ (cf. \cite[Theorem 11.7]{Ru} for harmonic functions).

We consider the associated Dirichlet boundary value problem:
\be\label{eq:492}
\begin{cases}
\begin{aligned}
L_{\alpha,\beta}u=0
&\text{ $~\mbox{on}~\mathbb{D};$}\\
u=f
&\text{ $~\mbox{on}~\mathbb{T}. $}
\end{aligned}
\end{cases}
\ee
Here the boundary function $f \in L^{1}(\mathbb{T})$,
and the boundary condition in \eqref{eq:492} is understood as $u_{r}\rightarrow f$ a.e. as $r\rightarrow 1^{-}$, where
$u_{r}(e^{i\theta})=u(re^{i\theta})$
for $e^{i\theta}\in \mathbb{T}$ and $r \in [0,1)$.
Klintborg and  Olofsson showed that  a function $u$ on $\mathbb D$ satisfies \eqref{eq:492} if and only if it has the form $u(z)=P_{\alpha,\beta}[f](z)$ for the boundary function $f \in C(\mathbb{T})$ (cf. \cite[Theorem 7.1]{Markus}).
We also mention that the Dirichlet problem for standard weighted Laplace differential operators $L_{0,\alpha}$ for arbitrary distributional boundary data was solved by Olofsson and Wittsten \cite{OLO}.
Subsequently, Carlsson and Wittsten \cite{Carlsson} solve the corresponding Dirichlet problem on the upper half plane by means of a counterpart of the classical Poisson integral formula.
For other related work, we refer the reader to \cite{Olo} due to Olofsson.


Before introducing the power series expansion of $(\alpha,\beta)$-harmonic functions, we recall the Gauss hypergeometric functions.
For $a,b,c\in \mathbb{C}$ such that $c\neq -1,-2,\ldots$,
the Gauss hypergeometric function is defined by the series
$$
F(a,b;c;x)=\sum_{n=0}^{\infty}\frac{(a)_n (b)_n}{(c)_n}\frac{x^n}{n\mbox{!}}
$$
for $|x|<1$, and by continuation elsewhere. Here $(a)_0=1$ and $(a)_{n} = a(a + 1)\ldots(a + n - 1)$ for $n = 1, 2, \ldots$  are the Pochhammer symbols.
Recall the following observations (cf. \cite{And}):

1. If $ \Re(c-a-b) > 0 $, then
\be\label{eq:498}
\lim_{x\rightarrow 1}F(a,b;c;x)=\frac{\Gamma(c)\Gamma(c-a-b)}{\Gamma(c-a)\Gamma(c-b)}.
\ee

2. It holds that
\be\label{eq:4204}
\frac { d F ( a , b ; c ; x ) } { d x } = \frac { a b } { c } F ( a + 1 , b + 1 ; c + 1 ; x ).
\ee

\blem$($\cite[Lemma 1.2]{Olo}$)$\label{eq:l2.6}
Let $c > 0,~a \leq c,~b \leq c$, and $ab \leq 0 ~(ab \geq 0)$. Then the function $F(a, b; c; x)$ is decreasing (increasing) on $(0, 1)$.
\elem

The Beta function, for real numbers $x>0$ and $y>0$, is defined as
$$B(x,y)=\int_{0}^{1}t^{x-1} (1-t)^{y-1}dt.$$
It has several useful properties such as $B(x,y)=B(y,x)$, and the relationship with the Gamma function $B(x,y)=\frac{\Gamma(x)\Gamma(y)}{\Gamma(x+y)}$.

\medskip

From \cite[Theorem 5.1]{Markus}, we know that a function $u$ on $\mathbb{D}$ is ($\alpha,\beta$)-harmonic if and only if it has the following convergent power series expansion:
\be\label{eq1.7}
u(z)=\sum_{k=0}^{\infty}c_{k}F(-\alpha,k-\beta;k+1;|z|^{2})z^{k}+
\sum_{k=1}^{\infty}c_{-k}F(-\beta,k-\alpha;k+1;|z|^{2})\overline{z}^{k},
\ee
where $\{c_k\}_{k=-\infty}^{\infty}$ denotes a sequence of complex numbers with
$$
\lim_{|k|\rightarrow\infty}\sup|c_k|^{1/|k|}\leq 1,
$$
and $F$ represents hypergeometric functions. 

Denote by $\mathcal{S}_{H}$ the class of all sense-preserving univalent harmonic functions $f=h+\overline{g}$ on $\mathbb{D}$ with the normalizations $h(0)=g(0)=h'(0)-1=0$. The class $\mathcal{S}_{H}^{0}:=\{f=h+\overline{g} \in \mathcal{S}_{H}:g'(0)=0\}$ is compact and in a one-to-one correspondence with $\mathcal{S}_{H}$.
The research on coefficient estimates of harmonic functions can be traced back to 1984, Clunie and Sheil-Small \cite{Clunie} gave the coefficient estimates for $f=h+\overline{g} \in \mathcal{S}_{H}^{0}$, where $h(z)=z+\sum_{n=2}^{\infty}a_{n}z^{n}$ and $g(z)=\sum_{n=2}^{\infty}b_{n}z^{n}$. More precisely, Clunie and Sheil-Small proved that $|a_2|<12172$ and proposed the following conjecture for all $n\geq 2$:
$$
|a_{n}|\leq \frac{(n+1)(2n+1)}{6} ~~\mbox{and}~~ |b_{n}|\leq \frac{(n-1)(2n-1)}{6}.
$$
Later, the inequality $|a_2|<49$ was proved in \cite{Duren}. Recently, Abu Muhanna et al. in \cite {Aub} further improved this estimate, and obtained $|a_2|<20.9197$. Furthermore, many scholars have studied on the coefficient estimates of different subclasses of harmonic functions. For example, the coefficient estimate for the special case of harmonic starlike functions have been confirmed, see \cite[P$_{107}$]{Duren} or \cite{Sh}. For convex harmonic functions, a sharper estimate is available (cf. \cite[P$_{50}$]{Duren} or \cite{Clunie}). For a coefficient estimate of $(\frac{\alpha}{2}, \frac{\alpha}{2})$-harmonic functions, we refer to \cite{BO2}. However, many fundamental questions concerning the coefficient estimates of harmonic functions remain open. This paper provides several estimates for the coefficients of ($\alpha,\beta$)-harmonic functions.

 The Koebe one-quarter theorem states that for every univalent analytic function $f$ defined on the unit disk, with $f(0)=0$ and $|f'(0)|=1$, its range contains the open disk $|w|<1/4$ (see \cite[P$_{32}$]{DUREN}). For a univalent sense-preserving harmonic function $f=h+\overline{g}$ on $\mathbb D$ satisfying the normalization $h(0)=g(0)=h'(0)-1=g'(0)=0$, the range of $f$ contains the disk $|w|<1/16$ (see \cite[P$_{92}$]{Duren} or \cite{Clunie}). In \cite{ShaoLin}, Chen and Ponnusamy established an asymptotically sharp Koebe type covering theorems for K-quasiconformal harmonic functions.
Little seems to be known about the covering theorems of $(\alpha,\beta)$-harmonic functions.
One aim of this paper is to establish a covering theorem to $(\alpha,\beta)$-harmonic functions.

In \cite{DUREN}, Duren proposed the growth and distortion estimates for univalent analytic functions. Later Clunie and Sheil-Small \cite{Clunie} studied the growth and distortion estimates for harmonic functions.  We refer to \cite{Khalfallah2} for the further improvement of these results on harmonic functions. In \cite{David}, Kalaj obtained some sharp or asymptotically sharp growth and distortion estimates for the real kernel $\alpha$-harmonic functions. There are also similar results for $\alpha$-harmonic functions in \cite{Jiao} and so-called complex-valued $\alpha$-harmonic functions in \cite{Khalfallah}.
Recently, Long \cite{Yong} gave  estimates of  real kernel $\alpha$-harmonic functions and their first order partial derivative functions by using the $L^p$ norm of the boundary functions.  For an $(\alpha,\beta)$-harmonic function $u$, Arsenovi\'{c} and Gaji\'{c} \cite{Milos} proved a sharp estimate of $|Du(0)|$ in terms of $L^p$ norm of the boundary function and gave the asymptotically sharp estimate of $|Du(z)|$.
 In \cite{PEI}, the authors considered the boundedness of $(\alpha,\beta)$-harmonic
functions under the condition that the boundary functions are bounded.
 The final aim of this paper is to generalize these results to the case of
 ($\alpha,\beta$)-harmonic functions.

The paper is organized as follows. In Section \ref{sec2}, we give coefficients estimates, Rad\'{o}'s theorem, Koebe type covering theorems, and an area theorem of $(\alpha,\beta)$-harmonic functions. 
In Section \ref{sec3}, we establish growth estimates and distortion estimates to $(\alpha,\beta)$-harmonic functions. 

\section{Coefficients estimates, Rad\'{o}'s theorem, Koebe type theorems and area theorem}\label{sec2}

\subsection{Coefficients estimates}

Heinz's lemma guarantees that the inequality
$$
|a_1|^2+\frac{3\sqrt{3}}{\pi}|a_0|^2+|b_1|^2\geq \frac{27}{4\pi^2}
$$
holds for the coefficients of any univalent harmonic mapping from the unit disk into itself.
The lower bound $27/4\pi^2$ is sharp (cf. \cite[p.\,67]{Duren}).
In \cite{BO2}, Long and Wang investigated Heinz-type inequalities for real kernel $\alpha$-harmonic functions:
$$
u_\alpha(z)=\sum_{k=0}^{\infty}c_{k}F\big (-\frac{\alpha}{2},k-\frac{\alpha}{2};k+1;|z|^{2}\big)z^{k}+
\sum_{k=1}^{\infty}c_{-k} F\big(-\frac{\alpha}{2},k-\frac{\alpha}{2};k+1;|z|^{2}\big )\overline{z}^{k}.
$$
They showed that
$$\left(\frac{\Gamma(\alpha+1)}{\Gamma(\frac{\alpha}{2}+1)^{2}}\right)^{2}\left(\frac{|c_{1}|^{2}}{(\frac{\alpha}{2}+1)^{2}}+\frac{3\sqrt{3}}{\pi}|c_{0}|^{2}
+\frac{|c_{-1}|^{2}}{(\frac{\alpha}{2}+1)^{2}}\right)\geq\frac{27}{4\pi^{2}},
$$
and the lower bound $\frac{27}{4\pi^{2}}$ is sharp.
For $\alpha\in(-1,0]$, if $u_\alpha$ is univalent and has real coefficients,
then the coefficient inequalities
$$
\left|\frac{c_{k}-c_{-k}}{1-c_{-1}}\right|\leq
\frac{\Gamma(k+1+\frac{\alpha}{2})}{\Gamma(k)\Gamma(2+\frac{\alpha}{2})}
$$
hold for $k=2,3,4,\ldots$. If $\alpha=0$, then the above inequalities reduce to the case of harmonic functions (cf. \cite[Theorem 6.4]{Clunie})

Now, we give Heinz's inequality of $(\alpha,\beta)$-harmonic  functions.

\bthm\label{thm2.02} Suppose that $\alpha$ and $\beta$ are non-negative integers.
Let $u$ be a univalent sense-preserving $(\alpha,\beta)$-harmonic function with the form \eqref{eq1.7}.
 Suppose that $u$ maps the unit disk $\mathbb{D}$ onto itself. Then,
\be\label{eq2.2}
\left(\frac{\Gamma(1+\alpha+\beta)}{\Gamma(1+\alpha)\Gamma(1+\beta)}\right)^2\left(\frac{|c_{1}|^2}{(1+\alpha)^2}
+\frac{3\sqrt{3}}{\pi}|c_{0}|^2+\frac{|c_{-1}|^2}{(1+\beta)^2}\right)\geq \frac{27}{4\pi^2}.
\ee
The lower bound is asymptotically sharp as $\alpha, \beta\rightarrow 0$.
\ethm

\bpf Suppose $u$ is of form \eqref{eq1.7}, univalent and sense-preserving on $\mathbb D$. By \eqref{eq1.7} and \eqref{eq:498},  for $z=re^{i\theta}\in \mathbb D$, we have
\beqq
e^{i\theta(t)}&:=&\lim_{r\rightarrow 1^{-}}u(z)\\
&=& \sum_{k=0}^{\infty}c_{k}F(-\alpha,k-\beta;k+1;1)e^{ikt}+
\sum_{k=1}^{\infty}c_{-k}F(-\beta,k-\alpha;k+1;1)e^{-ikt}\\
&=&\sum_{k=0}^{\infty}c_{k}\frac{\Gamma(1+\alpha+\beta)\Gamma(k+1)}{\Gamma(1+\beta)
\Gamma(k+1+\alpha)}e^{ikt}+
\sum_{k=1}^{\infty}c_{-k}\frac{\Gamma(1+\alpha+\beta)\Gamma(k+1)}{\Gamma(1+\alpha)
\Gamma(k+1+\beta)}e^{-ikt},
\eeqq
where $\theta(t)$ is a continuous nondecreasing function with
$\theta(t+2\pi)=\theta(t)+2\pi$.
By using the similar argument as that in  Hall's proof of Heinz's inequality (cf. \cite{Ha}),
we obtain the required inequality.
\epf

We note that Theorem \ref{thm2.02} generalizes the corresponding result in \cite{BO2}.

\bigskip

Rad\'{o} and Kneser in 1926 proved the following general result (cf. \cite{Kneser}):

\blem\label{lem2.1}
Let $\Omega$ be a bounded simply connected Jordan domain. Consider a homeomorphism $g^{*}$  from $\partial\mathbb{D}$ onto $\partial \Omega$. Let $g$ be the harmonic extension of $g^{*}$ into the disk. If $g(\mathbb{D})\subset \Omega$, then $g$ is univalent on $\mathbb{D}$.
\elem

Motivated by \cite[Lemma 1]{Daoud} for  polyharmonic functions, we give the
the following result, which is very important in the proof of our main results: 

\blem\label{lem2.2}
A locally univalent sense-preserving $(\alpha,\beta)$-harmonic function on the unit disk $\mathbb{D}$ of the form
$$
u(z)=\sum_{k=0}^{\infty}c_{k}F(-\alpha,k-\beta;k+1;|z|^{2})z^{k}+
\sum_{k=1}^{\infty}c_{-k}F(-\beta,k-\alpha;k+1;|z|^{2})\overline{z}^{k}
$$
is one-to-one if and only if, for each $r\in (0,1)$, the harmonic function
\be\label{eq:6181}
g_{r}(z)=\sum_{k=0}^{\infty}c_{k}F(-\alpha,k-\beta;k+1;r^{2})z^{k}+
\sum_{k=1}^{\infty}c_{-k}F(-\beta,k-\alpha;k+1;r^{2})\overline{z}^{k}
\ee
is a univalent sense-preserving harmonic function on $\mathbb{D}_{r}$.
\elem

\bpf
Assume that $u$ is univalent on $\mathbb{D}$ and set $u_{r}(z)=u(z)|_{\mathbb{D}_{r}}$. An initial observation reveals that
 for $r\in (0,1)$, $u_{r}(\mathbb{D}_{r})$ is a chain of outward open domains $\Omega_{r}$, whose boundary functions $f_{r}$ are univalent and disconnected, i.e., for $0<r<R\leq 1$,
$$
f_{r}(\partial \mathbb{D}_{r})\cap f_{R}(\partial \mathbb{D}_{R})=\emptyset .
$$
Moreover,
\be\label{eq:611a1}
u_{r}(\partial \mathbb{D}_{r})\equiv g_{r}(\partial \mathbb{D}_{r}).
\ee

Assume to the contrary that for some $0<r<1$, $g_{r}$ is not univalent on $\mathbb{D}_{r}$. Because $g_{r}$ is a univalent harmonic function on $\partial \mathbb{D}_{r}$, by Lemma \ref{lem2.1}, at some interior point, $\zeta\in\mathbb{D}_{r}$, $g_{r}(\zeta)\notin u_{r}(\mathbb{D}_{r})$. As such, $|\zeta|<r$ and there exists $\rho >r$ such that $g_{r}(\zeta) \in u_{\rho}(\partial \mathbb{D}_{\rho})$. Accordingly, $r>|\zeta|=\rho$, which contradicts $\rho >r$.

Now we assume to the contrary that for some $0<\rho<1$, $u_{\rho}$ is not univalent. Since $u_{\rho}$ is sense-preserving, there exists a point $\omega_0$ such that the total change of the argument of $u_{\rho}(z)-\omega_0$ around $u_{\rho}(\partial \mathbb{D}_{\rho})$ is
$$
\frac{1}{2\pi}\int_{0}^{2\pi} \arg(u_{\rho}(\rho e^{it})-\omega_0)\,dt \geq \,2.
$$
By the inclusion property, $u_{\rho}(\mathbb{D}_{\rho})\subset u(\mathbb{D})$, the total change of the argument of $u(z)-\omega_0$ around $u(\mathbb{T})$ is greater or equal to 2. Each $g_{r}$, $0<r<1$, is sense-preserving and univalent on $\mathbb{D}_{r}$. A limiting process of \eqref{eq:611a1} as $r\rightarrow 1$ shows that $g(z)\equiv g_{1}(z)$ is univalent on $\mathbb{D}$ and $g(z)\equiv u(z)$ for all $z\in \mathbb{T}$. Hence,
$$
1=\frac{1}{2\pi}\int_{0}^{2\pi} \arg(g(e^{it})-\omega_0)\,dt =\frac{1}{2\pi}\int_{0}^{2\pi} \arg(u(e^{it})-\omega_0)\,dt\geq 2,
$$
which leads to a contradiction.
\epf

By Lemma \ref{lem2.2}, the univalence of an $(\alpha,\beta)$-harmonic function $u$ on $\mathbb D$ is equivalent to the univalence of the harmonic functions $g_{r}$ on $\mathbb{D}_{r}$ with $0<r<1$. Clearly,
$$
\begin{aligned}
&g_{r}(rz)=\sum_{k=0}^{\infty}c_{k}F(-\alpha,k-\beta;k+1;r^{2})r^{k}z^{k}+
\sum_{k=1}^{\infty}c_{-k}F(-\beta,k-\alpha;k+1;r^{2})r^{k}\overline{z}^{k}\\
=&c_0F(-\alpha,-\beta;1;r^{2})+c_1F(-\alpha,1-\beta;2;r^{2})rz+\sum_{k=2}^{\infty}c_{k}F(-\alpha,k-\beta;k+1;r^{2})r^{k}z^{k}\\
&+c_{-1}F(-\beta,1-\alpha;2;r^{2})r\overline{z}+\sum_{k=2}^{\infty}c_{-k}F(-\beta,k-\alpha;k+1;r^{2})r^{k}\overline{z}^{k}.
\end{aligned}
$$
Let $A_k(r)=c_{k}F(-\alpha,k-\beta;k+1;r^{2})r^{k}$ and $B_k(r)=c_{-k}F(-\beta,k-\alpha;k+1;r^{2})r^{k}$.
For a sense-preserving $(\alpha,\beta)$-harmonic function $u$ with $u(0)=0$, it is easy to verify that $A_1(r)\neq 0$ for each $r\in [0,1)$. We let
\be\label{eq:6102}
F_{r}(z)=\frac{g_{r}(rz)}{A_1(r)}=z+\sum_{k=2}^{\infty}\frac{A_{k}(r)}{A_1(r)}z^{k}
+\sum_{k=1}^{\infty}\frac{B_{k}(r)}{A_1(r)}\overline{z}^{k}.
\ee

Denote by $\mathcal{S}_{(\alpha,\beta)}$ the class of all univalent sense-preserving $(\alpha,\beta)$-harmonic functions on $\mathbb D$ with the normalizations $c_0=c_{1}-1=0$. The subclass of functions $u \in \mathcal{S}_{(\alpha,\beta)}$ satisfying the additional condition $c_{-1} = 0$ is denoted by $\mathcal{S}_{(\alpha,\beta)}^0$.

Similar to the study of typically real functions in the theory of analytic functions or harmonic functions (cf. \cite[Section 6.6]{Duren}),
another interesting problem is the study of $(\alpha,\beta)$-harmonic functions with real coefficients.
Let $T\mathcal{S}_{(\alpha,\beta)}$ be the subclass of $\mathcal{S}_{(\alpha,\beta)}$ consisting of functions $u$ of the form \eqref{eq1.7} with  real coefficients. A function in $T \mathcal{S}_{(\alpha,\beta)}$ is said to be a typically real $(\alpha,\beta)$-harmonic function.

\begin{thm}\label{thm2.4}Suppose that $\alpha$ and $\beta$ are not negative integers.
Let $u\in T\mathcal{S}_{(\alpha,\beta)}$. Then
$$
\left|
\frac{\Gamma(1+\alpha)c_k}{\Gamma(k+1+\alpha)}- \frac{\Gamma(1+\beta)c_{-k}}{\Gamma(k+1+\beta)}\right|\leq
\frac{1}{\Gamma(k)}\left|\frac{1}{1+\alpha}-\frac{c_{-1}}{1+\beta}\right|
$$
for $k=2,3,4,\ldots$.
\end{thm}

\bpf Since $u\in \mathcal{S}_{(\alpha,\beta)}$, it follows from Lemma \ref{lem2.2} that, for each $r\in (0,1)$,  $g_r$ is univalent
and sense-preserving on $\mathbb D_r$. Obviously, $g_r$ has real coefficients, and then it is easy to verify that $g_r$
is typically real. Therefore,
\beq
\nonumber &&H_r(z):=\big(F(-\alpha,1-\beta;2;r^{2})-c_{-1}F(-\beta,1-\alpha;2;r^{2})\big)rz\\
&&
\nonumber \hspace{1.5cm}+\sum_{k=2}^{\infty}\big(c_{k}F(-\alpha,k-\beta;k+1;r^{2})-c_{-k}F(-\beta,k-\alpha;k+1;r^{2})\big)r^{k}z^{k}
\eeq
is analytic and typically real, since  $\Im{g_r}=\Im{H_r}$. By using \cite[p.\,58, Corollary]{DUREN}, we obtain, for $k\geq2$,
$$
\left|\frac{\big(c_{k}F(-\alpha,k-\beta;k+1;r^{2})-c_{-k}F(-\beta,k-\alpha;k+1;r^{2})\big)r^{k}}{\big(F(-\alpha,1-\beta;2;r^{2})
-c_{-1}F(-\beta,1-\alpha;2;r^{2})\big)r}\right|\leq k.
$$
This implies
$$
\left|
\frac{\Gamma(1+\alpha)c_k}{\Gamma(k+1+\alpha)}- \frac{\Gamma(1+\beta)c_{-k}}{\Gamma(k+1+\beta)}\right|\leq
\frac{1}{\Gamma(k)}\left|\frac{1}{1+\alpha}-\frac{c_{-1}}{1+\beta}\right|,
$$
which is the required inequality.
The proof is complete.
\epf

Theorem \ref{thm2.4} is a generalization of \cite[Theorem 1.3]{BO2}.

\bigskip

\blem$($\cite[Theorem 4.4]{Heikkala}$)$\label{lem2.3}
Let $r_{n}$ and $s_{n}~ (n = 0, 1, 2,...)$ be real numbers, and let the power series
$$ R ( x ) = \sum _ { n = 0 } ^ { \infty } r _ { n } x ^ { n }~~\mbox{and}~~S ( x ) = \sum _ { n = 0 } ^ { \infty } s _ { n } x ^ { n }$$
be convergent for $|x| < r$ $(r > 0)$ with $s_{n}> 0$ for all $n$. If the non-constant sequence $\{r_{n}/s_{n}\}$ is increasing (resp. decreasing) for all $n$, then the function $x\rightarrow R(x)/S(x)$ is strictly increasing (resp. decreasing) on $(0, r)$.
\elem

\blem\label{lem2.4}
For the integer $k\geq 1$ and the real numbers $\alpha$ and $\beta$ such that $\alpha+\beta>-1$, let $F_{k}(t)= F(-\alpha, k-\beta; k + 1 ; t )$ and $E_{k}(t)= F(-\beta, k-\alpha; k + 1 ; t )$ with $t\in (0,1)$.

$(1)$ If $\alpha=0$ or $k = 1$, then $\frac {F_{k}(t)} {F_{1}(t)} \equiv 1$ for $t \in(0, 1)$; If  $\alpha<0$ is not an integer, and $\beta\in(-1,1)$, then $\frac{ F _ { k }(t) } { F _ { 1 } (t)}$ is strictly increasing for $t \in(0, 1)$.

$(2)$ If $\alpha=0$ and $-1<\beta\leq 0$ , then $ \frac{E_{k}(t)}{F_{1}(t)} = E_{k}(t)$ is increasing for $t \in(0, 1)$;
 If $\alpha=0$ and $\beta\geq 0$ , then $ \frac{E_{k}(t)}{F_{1}(t)} = E_{k}(t)$ is decreasing for $t \in(0, 1)$;
 If $-1<\beta<\alpha<0$, then $\frac{E_{k}(t)}{F _{1}(t)}$ is strictly increasing for $t \in(0, 1)$.
\elem

\begin{proof}
(1) If $\alpha= 0$ or $k = 1$, then $F_1 (t)= F_{k} (t)\equiv1$.

Suppose $\alpha<0,~\beta<1$. Let
$$
A_{n }=\frac {(-\alpha)_{n}(k-\beta)_{n}}{(k+1) _ { n } n ! }~~\mbox{and}~~B_{n }=\frac {(-\alpha)_{n}(1-\beta) _ { n } } { (2) _ { n } n ! }
$$
for $n=0, 1, 2, \ldots $. Then it follows that $B_{n} > 0$ and
$$
\frac{A_{n}} { B_{n}} = \frac{(k-\beta)_{n} (2)_{n}}{(k+1)_{ n } (1-\beta)_{n}}.
$$
Since $\beta>-1$, it holds that
$$
\frac {{ A _ { n + 1}}/{B_{ n + 1}}}{{A_{n}}/{B_{n}}}= \frac {(k-\beta+ n ) (2+ n )} {( k+1+ n) (1-\beta+ n ) } > 1 .
$$
Thus ${A _ { n } } /{ B _ { n } }$ is strictly increasing for all $n$. By Lemma \ref{lem2.3}, we see that
$$
\frac{F _ { k } (t)} { F _ { 1 }(t) }=\frac{\sum _ { n = 0 }^{ \infty } A _ { n } t ^ { n } } { \sum _ { n = 0 } ^ { \infty } B _ { n } t ^ { n } }
$$
is strictly increasing for $t \in(0, 1)$.

(2) If $\alpha=0$, then $\frac{E_{k} (t)}{F_1 (t)}\equiv E_{k} (t)=F(-\beta,k;k+1;t)$. By applying Lemma \ref{eq:l2.6}, we obtain the desired result for the case $\alpha=0$.

Now, we assume that $\beta<\alpha<0$. Let
$$
C_{n }=\frac{(-\beta)_{n}(k-\alpha)_{n}}{(k+1)_{n}n!}
$$
for $n=0, 1, 2 \ldots $. Then it follows that $B_{n} > 0$ and
$$
\frac{C_n}{B_n}
= \frac{(k-\alpha)_n(-\beta)_n(2)_n}{(k+1)_n(-\alpha)_n(1-\beta)_n}.
$$
Since $-1<\beta<\alpha<0$, it holds that
$$
\frac { {C_ { n + 1}}/{B_{ n + 1}}}{{C_{n}}/{B_{n}}}= \frac{(k-\alpha+n)(-\beta+n)(2+n) }{(k+1+n)(-\alpha+n) (1-\beta+ n ) }  > 1 .
$$
Thus ${C_{n}}/{B_{n}}$ is strictly increasing for all $n$. By Lemma \ref{lem2.3}, we see that
$$
\frac{E_{k}(t)}{F_{1}(t)}=\frac{\sum_{n = 0}^{\infty} C_{n}t^{n}}{\sum_{n=0}^{\infty} B_{n} t^{n}}
$$
is strictly increasing for $t \in(0, 1)$.
\end{proof}

For the harmonic function $ f\in \mathcal S_{H}^{0}$ with the canonical representation $f=h+\overline{g}$, where $h(z)=z+\sum_{n=2}^{\infty}a_{n}z^{n}, g(z)=\sum_{n=2}^{\infty}b_{n}z^{n}$, Duren proved the sharp inequality $|b_2|\leq \frac{1}{2}$ (see \cite[Theorem, P$_{87}$]{Duren}), and Muhanna et al. \cite {Aub} obtained $|a_2|<20.9197$.
Next, we generalize the above results into the following form.

\bthm\label{thm-2.3}
Let $u$ belong to $\mathcal{S}_{(\alpha,\beta)}^{0}$, where $u$ has the form \eqref{eq1.7}, $-1<\beta<\alpha<0$. Then we have

$(1)$ $|c_{-2}|\leq \frac{(2+\beta)(1+\beta)}{4(1+\alpha)}$;

$(2)$ $|c_2|< 20.9197\left(1+\frac{\alpha}{2}\right)$.
\ethm

\begin{proof}
Based on Lemma \ref{lem2.2}, for $u \in \mathcal {S}_{(\alpha,\beta)-H}^{0}$, we obtain that the harmonic function $F_{r}\in \mathcal S_{H}^{0}$, where $F_{r}$ has the form \eqref{eq:6102}.

(1) According to \cite[Theorem, P$_{87}$]{Duren}, we obtain that
$$
\left|\frac{B_2(r)}{A_1(r)}\right|=\left|\frac{c_{-2}F(-\beta,2-\alpha;3;r^{2})r^{2}}{F(-\alpha,1-\beta;2;r^{2})r}\right|\leq\frac{1}{2}.
$$
It follows from Lemma \ref{lem2.4} that
$$
|c_{-2}|\leq \frac{1}{2}\inf_{r\in(0,1)}\frac{F(-\alpha,1-\beta;2;r^{2})}{F(-\beta,2-\alpha;3;r^{2})}= \frac{(2+\beta)(1+\beta)}{4(1+\alpha)}.
$$

If $\alpha=\beta=0$, then $|c_{-2}|\leq \frac{1}{2}$ holds, which shows that $\frac{(2+\beta)(1+\beta)}{4(1+\alpha)}$ is asymptotically optimal.

(2)
By using \cite[Theorem 1]{Aub}, we obtain that
$$
\left|\frac{A_2(r)}{A_1(r)}\right|=\left|\frac{c_{2}F(-\alpha,2-\beta;3;r^{2})r^{2}}{F(-\alpha,1-\beta;2;r^{2})r}\right|< 20.9197.
$$
Then  \eqref{eq:498} with Lemmas \ref{eq:l2.6} and \ref{lem2.4} yields that
$$
|c_2|<20.9197 \inf_{r\in (0,1)}\frac{F(-\alpha,1-\beta;2;r^{2})}{F(-\alpha,2-\beta;3;r^{2})} = 20.9197\left(1+\frac{\alpha}{2}\right).
$$
\end{proof}

\medskip

For the  harmonic functions  $f\in \mathcal S_{H}^{0}$, Clunie and Sheil-Small \cite{Clunie} proposed the following conjecture for all $n\geq 2$ :
$$
|a_{n}|\leq |A_{n}|=\frac{1}{6}(2n+1)(n+1)~~\mbox{and}~~|b_{n}|\leq |B_{n}|=\frac{1}{6}(2n-1)(n-1).
$$
If $u\in \mathcal{S}_{(\alpha,\beta)}^0$, which implies  $F_r\in S_{H}^{0}$, then we can suggest that
$$
\left|\frac{A_{k}(r)}{A_1(r)}\right|\leq \frac{1}{6}(2k+1)(k+1),~~\left|\frac{B_{k}(r)}{A_1(r)}\right|\leq \frac{1}{6}(2k-1)(k-1)
$$
for all indices $k \geq 2$. This implies
the coefficients conjecture for $(\alpha,\beta)$-harmonic functions.

\bcon\label{conj1}
Let $u\in \mathcal{S}_{(\alpha,\beta)}^{0}$ be a univalent $(\alpha,\beta)$-harmonic function, where $u$ has the form \eqref{eq1.7}. Then
$$
|c_k|\leq \frac{1}{6}(2k+1)(k+1)\inf_{r\in(0,1)} \frac{F(-\alpha,1-\beta;2;r^{2})}{F(-\alpha,k-\beta;k+1;r^{2})},$$$$
|c_{-k}|\leq \frac{1}{6}(2k-1)(k-1)\inf_{r\in(0,1)} \frac{F(-\alpha,1-\beta;2;r^{2})}{F(-\beta,k-\alpha;k+1;r^{2})}.
$$
\econ

A version of this  bound for $K$-quasiconformal harmonic mapping that is expected to be sharp
for all $K\geq 1$ has been discussed in \cite{Wa}.

\medskip

By using Lemma \ref{lem2.4}, if $-1<\beta<\alpha<0$,  we have
$$\inf_{r\in(0,1)} \frac{F(-\alpha,1-\beta;2;r^{2})}{F(-\alpha,k-\beta;k+1;r^{2})}
=\frac{\Gamma(k+1+\alpha)}{k!\Gamma(2+\alpha)}$$
and
$$\inf_{r\in(0,1)} \frac{F(-\alpha,1-\beta;2;r^{2})}{F(-\beta,k-\alpha;k+1;r^{2})}
=\frac{\Gamma(k+1+\beta)}{(1+\alpha)k!\Gamma(1+\beta)}.$$

 For $\alpha$-harmonic functions $u\in \mathcal{S}_{(0,\alpha)-H}^{0}$ of the form \eqref{eq1.7},  we conjecture that
$$
|c_{k}|\leq \frac{1}{6}(2k+1)(k+1),$$
$$~~|c_{-k}|\leq \begin{cases}
\begin{array}{rl}
\frac{(2k-1)(k-1)}{6},& \ \alpha\geq 0;\\
\frac{(2k-1)(k-1)\Gamma(k+1+\alpha)}{6 k! \Gamma(1+\alpha)},& \  -1<\alpha<0.\\
\end{array}
\end{cases}
$$

By \cite[Theorem, P$_{107}$]{Duren} (or \cite{Sh}), we see that the coefficients of every starlike function $f \in \mathcal{S}_{H}^{0}$ satisfy the sharp inequalities:
\be\label{eq2.1}
|a_{n}|\leq\frac{1}{6}(2n+1)(n+1),~~~~|b_{n}|\leq\frac{1}{6}(2n-1)(n-1).
\ee
Based on the above inequalities, we can deduce the following results:

\bthm\label{thm2.2} Suppose that $\alpha,\beta \notin \{-1, -2\}$.
Let $u\in \mathcal{S}_{(\alpha,\beta)}^{0}$ be a starlike function, where $u$ is of form \eqref{eq1.7}. 
Then, for $k\geq 2$, we have
$$
|c_{k}|\leq \frac{(2k+1)(k+1)\Gamma(k+1+\alpha)}{6k!\Gamma(2+\alpha)},~|c_{-k}|\leq \frac{(2k-1)(k-1)\Gamma(k+1+\beta)}{6(1+\alpha)k!\Gamma(1+\beta)}.
$$
\ethm

\begin{proof}
It is easy to verify that $u(|z|=r)=g_r(|z|=r)$ and $u(\mathbb D_r)=g_r(\mathbb D_r)$,
where $g_r$ is defined by \eqref{eq:6181} in Lemma \ref{lem2.2}.
The radial limit function  $g_1(z)=\lim_{r\rightarrow 1}g_r(rz)$
is well defined. Since each $g_r(rz)$ is univalent and sense-preserving, by using
 the Hurwitz's theorem of harmonic functions (cf. \cite[p.\,10]{Duren}), we obtain that $g_1$ is a univalent sense-preserving harmonic function. In addition, $g_1(z)\equiv u(z)$ for all $z\in \partial\mathbb{D}$ and
 $u(\mathbb D)=g_1(\mathbb D)$, which shows
 $F_1\equiv \frac{g_1}{A_1(1)} \in \mathcal{S}_{H}^{0}$ is starlike.  By using \eqref{eq2.1},
$$
\left|\frac{A_{k}(1)}{A_1(1)}\right|\leq \frac{1}{6}(2k+1)(k+1),~~\left|\frac{B_{k}(1)}{A_1(1)}\right|\leq \frac{1}{6}(2k-1)(k-1).
$$
This implies the desired inequalities.

\end{proof}

\subsection{Rad\'{o}'s theorem}

There is no univalent harmonic function from $\mathbb{D}$ onto $\mathbb{C}$. This fact
 is the famous Rad\'{o}'s theorem (cf. \cite[P$_{24}$]{Duren}). Now, we give  Rad\'{o}'s theorem for $(\alpha,\beta)$-harmonic functions.
  A version of this result for polyharmonic mappings is given in \cite{Daoud}.

\bthm\label{thm-2.5}Let $\alpha$ and $\beta$ be real numbers such that
$\alpha+\beta>-1$ and $\alpha,\beta \notin \{-1, -2\}$. Then there is no univalent $(\alpha,\beta)$-harmonic function $u$ with $u(0)=0$ of the unit disk onto the whole complex plane.
\ethm

\begin{proof}
Suppose that a univalent $(\alpha,\beta)$-harmonic function $u$ with $u(0)=0$  maps $\mathbb{D}$ onto $u(\mathbb D)$  which contains a disk ${\mathbb D}_{R}=\{w \in \mathbb{C}: |w|< R\}$. 
By Lemma \ref{lem2.2}, we have $u(\mathbb{D}_{r})=g_{r}(\mathbb D_r)$ for $0<r<1$,
where $g_r$ is defined by \eqref{eq:6181}.
Thus, for arbitrary constant $\varepsilon\in (0,1)$,  we can assume that there
exists an $r_0=r_0(\varepsilon)$ such that, for each $r$ with $r_{0}\leq r<1$, $g_{r}(\mathbb{D}_{r})\supset {\mathbb D}_{(1-\varepsilon)R}$, that is $g_{r}(r \mathbb{D})\supset {\mathbb D}_{(1-\varepsilon)R}$.
For the harmonic function $g_{r}(rz)$, by using the similar arguments as that in the proof of \cite[Rad\'{o} Theorem, P$_{24}$]{Duren}, we can conclude that
$$c (1-\varepsilon)^2 R^2\leq |r(g_{r})_{z}(0)|^2+|r(g_{r})_{\overline{z}}(0)|^2,$$
where c is the Heinz constant $27/4\pi^{2}$.
It follows that
$$c(1-\varepsilon)^2 R^2\leq A_1(r)^2+B_1(r)^2.
$$
 Then, when $r\rightarrow 1$ and $\varepsilon\rightarrow 0$, we obtain
$$
cR^2 \leq \left(\frac{|c_1|\Gamma(1+\alpha+\beta)}{|\Gamma(2+\alpha)\Gamma(1+\beta)|}\right)^2
+\left(\frac{|c_{-1}|\Gamma(1+\alpha+\beta)}{|\Gamma(1+\alpha)\Gamma(2+\beta)|}\right)^2.
$$
Hence $R$ is  bounded. In particular, the range of $u$ cannot contain disks of arbitrarily large radius centered at the origin.
\end{proof}

\bcor[Rad\'{o}'s theorem for $\alpha$-harmonic functions]
There is no univalent $\alpha$-harmonic function of the unit disk onto the whole complex plane.
\ecor

\subsection{Koebe type theorems}

In the following, we consider Koebe type covering theorems for $(\alpha,\beta)$-harmonic functions. First, recall the corresponding ones for harmonic functions.

\begin{Thm}$($\cite[Theorem 1, P$_{90}$]{Duren}$)$ \label{thmA}
Each function in $\mathcal S_{H}$ omits some point on the circle $|w|=\frac{2\pi \sqrt{6}}{9}$. Each function in $\mathcal S_{H}^{0}$ omits some point on the circle $|w|=\frac{2\pi \sqrt{3}}{9}$, but need not omit any point of smaller modulus.
\end{Thm}

\begin{Thm}$($\cite[Theorem 4.4]{Clunie}$)$\label{thmB}
Each function $f \in \mathcal S_{H}^{0}$ satisfies the inequality
$$
|f(z)|\geq \frac{1}{4}\frac{|z|}{(1+|z|)^{2}},~~|z|<1.
$$
In particular, the range of $f$ contains the disk $|w|<\frac{1}{16}$.
\end{Thm}

\bthm\label{thm-2.6} Suppose that $\alpha,\beta \notin \{-1, -2\}$.
 If $u$ belongs to $\mathcal{S}_{(\alpha,\beta)}$, then $u$ omits some point on the circle $$|w|=\frac{2\pi \sqrt{6}}{9}\frac{\Gamma(1+\alpha+\beta)}{|\Gamma(2+\alpha)\Gamma(1+\beta)|}.$$
If $u\in \mathcal{S}_{(\alpha,\beta)}^{0}$, then $u$ omits some point on the circle $$|w|=\frac{2\pi \sqrt{3}}{9}\frac{\Gamma(1+\alpha+\beta)}{|\Gamma(2+\alpha)\Gamma(1+\beta)|}.$$
\ethm

\begin{proof}
 In the proof of Theorem \ref{thm-2.5}, it was shown that if an $(\alpha,\beta)$-harmonic function $u$ with $u(0)=0$ contains a disk $|w|<R$ in its range, then there exists $r_0\in (0,1)$ such that, for each $r$ with $r_0<r<1$,
$$
c(1-\varepsilon)^2R^2\leq A_1(r)^2+B_1(r)^2,
$$
where  $c=27/4\pi^{2}$. Thus,
$$
R^2\leq \frac{4\pi^{2}}{27(1-\varepsilon)^2}\big(A_1(r)^2+B_1(r)^2\big).
$$
If $u \in \mathcal{S}_{(\alpha,\beta)}$, then $c_1=1$ and $|B_1(r)|\leq |A_1(r)|$, so when $r\rightarrow 1$ and $\varepsilon\rightarrow 0$, it implies that $$R\leq\frac{2\sqrt{6}\pi}{9}\frac{\Gamma(1+\alpha+\beta)}{|\Gamma(2+\alpha)\Gamma(1+\beta)|}.$$
If $u \in \mathcal{S}_{(\alpha,\beta)}^{0}$, then $c_1=1$ and $c_{-1}=0$, so when $r\rightarrow 1$ and $\varepsilon\rightarrow 0$, it follows that $$R\leq\frac{2\sqrt{3}\pi}{9}\frac{\Gamma(1+\alpha+\beta)}{|\Gamma(2+\alpha)\Gamma(1+\beta)|}.$$ In either case it follows that some omitted value must lie on the circle of the given radius.

\end{proof}

For $\alpha=\beta=0$, we have  Theorem A for harmonic functions which shows Theorem \ref{thm-2.6} is asymptotically sharp as $\alpha, \beta\rightarrow 0$  for the class $\mathcal{S}_{(\alpha,\beta)}^{0}$.

\bthm\label{thm-2.7}
Let $u$ belong to $\mathcal{S}_{(\alpha,\beta)}^{0}$ with $\alpha,\beta \notin \{-1, -2\}$. Then $u$ satisfies the inequality
$$
|u(z)|\geq \frac{|A_1(|z|)|}{16},~~|z|<1.
$$
In particular, the range of $u$ contains the disk $|w|<\frac{1}{16}\frac{\Gamma(1+\alpha+\beta)}{|\Gamma(2+\alpha)\Gamma(1+\beta)|}$.
\ethm

\begin{proof}
Since $u \in \mathcal{S}_{(\alpha,\beta)}^{0}$, by Lemma \ref{lem2.2}, we have $F_{r} \in \mathcal S_{H}^{0}$, where $F_{r}$ has the form \eqref{eq:6102}. By Theorem B, we have
$$
|F_{r}(z)|=\left|\frac{g_{r}(rz)}{A_1(r)}\right|\geq \frac{1}{4}\frac{|z|}{(1+|z|)^{2}}.
$$
Hence, for $|z|=r$,
$$
|u(z)|=|g_{r}(z)| \geq 
\frac{|A_1(r)|}{16}.
$$
Moreover,
$$
\lim_{|z|=r\rightarrow 1}|u(z)|\geq \frac{|A_1(1)|}{16}=\frac{1}{16}\frac{\Gamma(1+\alpha+\beta)}{|\Gamma(2+\alpha) \Gamma(1+\beta)|}.
$$
\end{proof}

\bcor
For an $\alpha$-harmonic function $u$ which belongs to $\mathcal{S}_{(0,\alpha)-H}^{0}$, we have
$$
|u(z)|\geq \frac{|z|}{16}.
$$
\ecor

\subsection{Area theorem}

Among all functions in $\mathcal S_{H}^{0}$, which ones map the unit disk to a region of smallest area? The following theorem gives an affirmative answer to this questions.

\begin{Thm}\label{thmC}
The area of the image of the unit disk under each function $f$ in $\mathcal S_{H}^{0}$ is greater than or equal to $\frac{\pi}{2}$, and the minimum is attained only by the function $f(z)=z+\frac{1}{2}\overline{z}^{2}$ and its rotations.
\end{Thm}

We improve and generalize Theorem C into the following form.

\bthm\label{thm-2.8} Let $u \in \mathcal{S}_{(\alpha,\beta)}^{0}$ with $\alpha,\beta \notin \{-1, -2\}$.
The area of the image of the unit disk under $u$ is greater than or equal to
 $$\frac{\pi }{2}\frac{\Gamma(1+\alpha+\beta)}{|\Gamma(2+\alpha)\Gamma(1+\beta)|},$$  which is asymptotically sharp as $\alpha, \beta\rightarrow 0$.
\ethm

\begin{proof} Suppose that $u \in \mathcal{S}_{(\alpha,\beta)}^{0}$ is of the form \eqref{eq1.7} with $c_0=c_{-1}=c_1-1=0$, and $g_r$ is defined as in \eqref{eq:6181}. By Lemma \ref{lem2.2}, for each $r\in (0,1)$, $g_r$ is univalent and sense-preserving.
The Jacobian of $g_{r}$ is $|(g_{r})_{z}|^{2}-|(g_{r})_{\overline{z}}|^{2}$, so the area of $g_{r}(\mathbb{D}_{r})$ is

$$
\begin{aligned}
\mathcal{A}(g_{r}(\mathbb{D}_{r}))&=\iint_{\mathbb{D}_{r}}\left(|(g_{r})_{z}(z)|^{2}-|(g_{r})_{\overline{z}}(z)|^{2}\right)\, dx dy\\
&\geq \iint_{\mathbb{D}_{r}}(1-|z|^{2})|(g_{r})_{z}(z)|^{2}\,dx dy,
\end{aligned}
$$
where $$(g_{r})_{z}(z)= \sum_{k=1}^{\infty} kc_{k}F(-\alpha,k-\beta;k+1;r^{2})z^{k-1}:= \sum_{k=1}^{\infty}C_{k}z^{k-1}.$$ 

Further,
$$
\mathcal{A}(g_{r}(\mathbb{D}_{r}))\geq 
\pi r^{2}(1-\frac{r^{2}}{2})|C_1|^2+\pi\sum_{k=2}^{\infty}(\frac{1}{k}-\frac{r^{2}}{k+1})|C_{k}|^{2}r^{2k}.
$$
Clearly, the last sum is minimized by choosing $C_{k}=0$ for all $k\geq 2$. Because $g_r$ and $u$ are both univalent, by using \eqref{eq:611a1}, $\ID_r$ has the same image under these mappings, and then we obtain
 $\mathcal{A}(u(\mathbb{D}_r))=\mathcal{A}(g_{r}(\mathbb{D}_{r}))$. Hence, when $r\rightarrow 1$, we have
$$
\mathcal{A}(u(\mathbb{D}))\geq \frac{\pi}{2}|A_1(1)|=\frac{\pi }{2}\cdot\frac{\Gamma(1+\alpha+\beta)}{|\Gamma(2+\alpha)\Gamma(1+\beta)|}.
$$
\end{proof}

\bcor
For each $\alpha$-harmonic function belong to $\mathcal{S}_{(0,\alpha)-H}^{0}$, the area of the image is greater than or equal to $\frac{\pi }{2(1+\alpha)}$, which is asymptotically sharp as $\alpha\rightarrow 0$.
\ecor

\section{Growth and distortion estimates}\label{sec3}

\subsection{Growth estimates}

In this subsection, we will consider the growth of $(\alpha,\beta)$-harmonic functions in terms of $L^p$ norm of the boundary functions. Now, we give the following notations.

Let $f$ be a measurable complex-valued function defined on $\mathbb{D}$. The integral means of $f$
is defined as follows:
$$
M _ {p} ( r , f ) = \left( \frac { 1 } { 2 \pi } \int _ { 0 } ^ { 2 \pi } | f ( r e ^ { i \theta } ) | ^ {
p } d \theta \right) ^ {\frac{1}{ p} } , ~~0 < p < \infty,
$$
and
$$
M _ { \infty } ( r , f ) = {\mbox{ess}\sup}_{ 0 \leq \theta \leq 2 \pi } | f ( r e ^ { i \theta } ) |,
$$
where $0<r<1$.
A function $f$ analytic on $\mathbb D$ is said to be of class $H^{p}(\mathbb{D})$ called the Hardy space if $M_{p}(r, f)$ is
bounded.
It is also convenient to define analogous classes $h^{p}(\mathbb{D})$ or $h^{p}_{\alpha,\beta}(\mathbb{D})$ of harmonic functions or  $(\alpha,\beta)$-harmonic functions, respectively. For more information about Hardy spaces, we refer to \cite{Du}.

Denote by $L^{p}(\mathbb{T})$, where $p \in [1,\infty]$, the space of all measurable functions $f$ of $\mathbb{T}$ into $\mathbb{C}$ with
$\|f\|_{L^{p}(\mathbb{T})}<\infty$, where
$$
\|f\|_{L^{p}(\mathbb{T})}=
\begin{cases}
\begin{array}{rl}
\left(\frac{1}{2\pi}\int_{0}^{2\pi}| f (e ^ { i \theta } ) | ^ {p} d \theta\right) ^ {\frac{ 1}{ p} },
&\text{ $p \in [1,\infty);$}\\
{\mbox{ess}\sup}_{\theta \in [0,2\pi]}| f (e ^ { i \theta } ) |,
&\text{$ p=\infty.$}
\end{array}
\end{cases}
$$

\blem $($\cite[P$_{1179}$]{Kalaj}$)$\label{lem3.2}
It holds that
\be\label{eq:6301}
\int _ { 0 } ^ { \pi } \frac { \sin ^ { \mu - 1 } t } { ( 1 + r ^ { 2 } - 2 r \cos t ) ^ { \nu } } d t
= B\Big( \frac { \mu } { 2 } , \frac { 1 } { 2 } \Big) F\Big(\nu , \nu + \frac { 1 - \mu } { 2 } ; \frac { 1 + \mu }{2}; r^{2} \Big),
\ee
where $B(u, v)$ is the Beta function, and $F(a, b; c; x)$ is the Gauss hypergeometric function. In particular,
$$\int _ { 0 } ^ { \pi } \sin ^ { \mu - 1 } t \left( 1 - \cos t \right) ^ { - \nu } d t = 2 ^ { \nu }
B \Big(\frac { \mu } { 2 } , \frac { 1 } { 2 } \Big) F\Big( \nu, \nu + \frac { 1 - \mu } { 2 }; \frac { 1 + \mu } { 2}; 1 \Big),
$$
and
\be\label{eq:4918}
\int_{0}^{\pi}\frac {dt} {(1+r^{2}-2r \cos t)^{\nu } } = \pi F ( \nu , \nu; 1; r ^ { 2 } ).
\ee
\elem

Next, we prove some asymptotically sharp results for the class of ($\alpha,\beta$)-harmonic functions.
Based on \cite[Theorem 1.1]{David} and \cite[Theorem 1.1]{Yong}, we arrive at the following  conclusion.

\bthm\label{thm3.1}
Let $u(z)=P_{\alpha,\beta}[f](z)$ be an $(\alpha,\beta)$-harmonic function on $\mathbb D$ s.t. the complex-valued function $f \in
L^{p}(\mathbb{T})$ $(p\geq1)$.
Then there is a function $A_{\alpha,\beta, p}(r)$  in $(0,1)$ and a constant $A_{\alpha,\beta, p}=\sup_{r\in(0,1)}
A_{\alpha,\beta, p}(r)$ defined in \eqref{eq:493} and \eqref{eq:494} below, so that
$$
|u(z)|\leq \frac{A_{\alpha,\beta, p}(r)}{(1-r^{2})^{1/p}} \|f\|_{L^{p}(\mathbb{T})} \leq \frac{A_{\alpha,\beta, p}}{(1-r^{2})^{1/p}} \|f\|_{L^{p}(\mathbb{T})},~~z\in\mathbb{D}.
$$
In particular, if $p=\infty$, we have $|u(z)|\leq \|f\|_{L^{\infty}(\mathbb{T})}$.

In addition, we have the following sharp inequality
\be\label{eq:4911}
M_{p}(r,u)\leq |c_{\alpha,\beta}|F\Big(-\frac{\alpha+\beta}{2},-\frac{\alpha+\beta}{2};1;r^{2}\Big)
\|f\|_{L^{p}(\mathbb{T})}
\ee
 with $r\in (0,1)$.
\ethm

\begin{proof}
We rewrite \eqref{eq:497} as the following  formula:
\be\label{eq:4910}
u(z)=P_{\alpha,\beta}[f](z)=\frac{1}{2 \pi}\int_{0}^{2\pi}P_{\alpha,\beta}(ze^{-it})f(e^{it})dt,
\ee
 where $z\in \mathbb D$.
Then, by the H\"{o}lder inequality, it follows that
$$
|u(z)|\leq\|f\|_{p}\left(\int_{0}^{2\pi}|P_{\alpha,\beta}(ze^{-it})|^{q}\frac{dt}{2\pi}\right)^{\frac{1}{q}}
$$
for $q$ such that $\frac{1}{p}+\frac{1}{q}=1$. We write $z=re^{i\theta}.$ For $t, s \in [0,2\pi]$, let
\be\label{eq:5281}
e ^ { i ( t - \theta ) } = \frac { r + e ^ { i s } } { 1 + r e ^ { i s } } .
\ee
It is easy to verify that
\be\label{eq:5282}
| 1 - r e ^ { i ( \theta - t ) } | = \frac { 1 - r ^ { 2 } } { | 1 + r e ^ { - i s } | }~\mbox{and}~d t = \frac { 1 - r ^ { 2 } } { | 1 + r e ^ { i s } | ^ { 2 } } d s.
\ee
Hence we obtain
$$
|P_{\alpha,\beta}(ze^{-it})|^{q}
=|c_{\alpha,\beta}|^{q}(1-r^{2})^{-q}(1+r^{2}+2r\cos s)^{\frac{(\alpha+\beta+2)q}{2}},
$$
and it follows that
$$
|u(z)|\leq |c_{\alpha,\beta}|\left(\int_{0}^{2\pi}A(q,s)\frac{ds}{2\pi}\right)^{\frac{1}{q}}\|f\|_{L^{p}(\mathbb{T})},
$$
where
$$
A(q,s)=(1-r^{2})^{1-q}(1+r^{2}+2r\cos s)^{\frac{(\alpha+\beta+2)q}{2}-1}.
$$
Now, let
$$
A_{\alpha,\beta,p}(r)=|c_{\alpha,\beta}|\left(\int_{0}^{2\pi}(1+r^{2}+2r\cos s)^{\frac{(\alpha+\beta+2)q}{2}-1}\frac{ds}{2\pi}\right)^{\frac{1}{q}}.
$$

Since $2-(2+\alpha+\beta)q<1$ holds for all $q\geq1$, and $\alpha+\beta>-1$,
by using \eqref{eq:4918}, Lemma \ref{eq:l2.6} and \eqref{eq:498} lead to
\be\label{eq:33001}
\begin{aligned}
&\int_{0}^{2\pi}(1+r^{2}+2r\cos s)^{\frac{(\alpha+\beta+2)q}{2}-1}\frac{ds}{2\pi}=
\frac{1}{\pi} \int _{ 0 } ^ { \pi } ( 1 + r ^ { 2 } - 2 r \cos t ) ^{\frac { (\alpha +\beta+2) q}{2}-1} dt\\
&= F\left(1-\frac{(\alpha+\beta+2)q}{2}, 1- \frac {(\alpha+\beta+2 ) q } {
2 } ; 1 ; r ^ { 2 } \right) \\
&<  F \left(1- \frac {(\alpha +\beta+2) q } { 2 } , 1-\frac {(
\alpha+\beta+2 ) q } { 2 } ; 1 ; 1 \right)\\
&  = \frac{\Gamma\big((\alpha+\beta+2)q-1\big)}{\Gamma^{2}
\Big(\frac{(\alpha+\beta+2)q}{2}\Big)}.
\end{aligned}
\ee
We obtain
\be\label{eq:493}
A_{\alpha,\beta, p}(r)=|c_{\alpha,\beta}|\left(
F\Big(1-\frac{(\alpha+\beta+2)q}{2}, 1-\frac {(\alpha+\beta+2 ) q } {2 } ; 1 ; r ^ { 2 } \Big)\right)^{\frac{1}{q}}
\ee
and
\be\label{eq:494}
A_{\alpha,\beta, p}=\sup_{r\in(0,1)}
A_{p,\alpha,\beta}(r)
=|c_{\alpha,\beta}|
\left(\frac{\Gamma\big((\alpha+\beta+2)q-1\big)}{\Gamma^{2}
\Big(\frac{(\alpha+\beta+2)q}{2}\Big)}\right)^{\frac{1}{q}}.
\ee
Therefore, we reduce
$$
|u(z)|\leq(1-|z|^{2})^{\frac{1}{q}-1}A_{\alpha,\beta, p}(r)\|f\|_{L^{p}(\mathbb{T})}\leq(1-|z|^{2})^{\frac{1}{q}-1}A_{\alpha,\beta, p}\|f\|_{L^{p}(\mathbb{T})}.
$$

Now, we prove the inequality \eqref{eq:4911}.  For
 $z=re^{i\theta}\in \mathbb D$, we have
\be\label{eq:499}
|u(z)|\leq \int_{0}^{2\pi} |c_{\alpha,\beta}|\frac{(1-r^{2})^{\alpha+\beta+1}}
{|1-re^{i(\theta-t)}|^{\alpha+\beta+2}}|f(e^{it})|\frac{dt}{2\pi}.
\ee
Let
$$
I= \int_{0}^{2\pi} |c_{\alpha,\beta}|\frac{(1-r^{2})^{\alpha+\beta+1}}
{|1-re^{i(\theta-t)}|^{\alpha+\beta+2}}\frac{dt}{2\pi}.
$$
Then from the proof of  \cite[Theorem 3.1]{Olo}, we obtain that
\beqq
&&I=|c_{\alpha,\beta}|F\Big(-\frac{\alpha+\beta}{2},-\frac{\alpha+\beta}{2};1;r^{2}\Big)
<|c_{\alpha,\beta}|F\Big(-\frac{\alpha+\beta}{2},-\frac{\alpha+\beta}{2};1;1\Big)\\
&&\hspace{0.35cm}=\frac{\Gamma(\alpha+1)\Gamma(\beta+1)}{\Gamma^{2}(1+\frac{\alpha+\beta}{2})},
\eeqq
where Lemma \ref{eq:l2.6} is used.

For $p\geq 1$, considering Jensen's inequality for \eqref{eq:499}, we have
$$
\begin{aligned}
|u(z)|^{p}&\leq\bigg(I\cdot \int_{0}^{2\pi}\frac{1}{I} \frac{|c_{\alpha,\beta}|(1-r^{2})^{\alpha+\beta+1}}
{|1-re^{i(\theta-t)}|^{\alpha+\beta+2}}|f(e^{it})|\frac{dt}{2\pi}\bigg)^{p}\\
&\leq I^{p-1}\int_{0}^{2\pi} |c_{\alpha,\beta}| \frac{(1-r^{2})^{\alpha+\beta+1}}
{|1-re^{i(\theta-t)}|^{\alpha+\beta+2}}|f(e^{it})|^{p}\frac{dt}{2\pi}.
\end{aligned}
$$
It follows that
$$
\begin{aligned}
\frac{1}{2\pi}\int_{0}^{2\pi}|u(z)|^{p}d\theta &\leq
I^{p-1}\frac{1}{2\pi}\int_{0}^{2\pi}\left(\int_{0}^{2\pi}|c_{\alpha,\beta}| \frac{(1-r^{2})^{\alpha+\beta+1}}
{|1-re^{i(\theta-t)}|^{\alpha+\beta+2}}|f(e^{it})|^{p}\frac{dt}{2\pi}\right)d\theta\\
& =I^{p}\|f\|_{L^{p}(\mathbb{T})}^{p}.
\end{aligned}
$$
Therefore,
$$
M_{p}(r,u)\leq I\|f\|_{L^{p}(\mathbb{T})}=|c_{\alpha,\beta}|F\Big(-\frac{\alpha+\beta}{2},-
\frac{\alpha+\beta}{2};1;r^{2}\Big)
\|f\|_{L^{p}(\mathbb{T})}.
$$

If $\alpha=\beta=0$, then $I\leq 1$, and the inequality $M_{p}(r,u)\leq\|f\|_{L^{p}(\mathbb{T})}$ holds, which gives \cite[Theorem 1.1]{Yong}. Let $f(e^{it})\equiv C$ with $C>0$. Then \eqref{eq:4910} shows that $u(z)=C c_{\alpha,\beta}F(-\frac{\alpha+\beta}{2},-\frac{\alpha+\beta}{2};1;r^{2})$. From this we see that inequality \eqref{eq:4911} is sharp.
\end{proof}

\subsection{Distortion estimates}

We begin this subsection with two lemmas.

\blem\label{lem3.3}
For $q\geq0$ and $r\in[0,1]$, let
$$
L(y)=\int_{0}^{2\pi}(A+B|\cos x|)^{q}(1+r^{2}+2r\cos(x-y))^{m}dx,
$$
where $y\in[0,2\pi]$, $A\geq 0$ and $B>0$ are constants. Then,
$$
\max_{y\in[0,2\pi]} L(y)=\max_{y\in[0,\pi]} L(y)=
\begin{cases}
\begin{array}{rl}
L(0),&
\text{$m>1;$}\\
L(\pi/2),&
\text{$ m\leq1.$}
\end{array}
\end{cases}
$$
\elem

\begin{proof}
Elementary computations show that
$$
\begin{aligned}
L^{'}(y)=&2mr\int_{0}^{2\pi}(A+B|\cos x|)^{q}(1+r^{2}+2r\cos(x-y))^{m-1}\sin(x-y)\, dx\\
=&2mr\int_{0}^{2\pi}(A+B|\cos (x+y)|)^{q}(1+r^{2}+2r\cos x)^{m-1}\sin x \, dx\\
=&2mr\int_{0}^{\pi}(A+B|\cos (x+y)|)^{q}\\
&\hspace{1cm}\times\big((1+r^{2}+2r\cos x)^{m-1}-(1+r^{2}-2r\cos x)^{m-1}\big)\sin x \,dx.
\end{aligned}
$$
Let
$$
C(x)=(1+r^{2}-2r\sin x)^{m-1}-(1+r^{2}+2r\sin x)^{m-1}.
$$
It follows that
 $$
\begin{aligned}
L^{'}(y)&=2mr\int_{-\pi/2}^{\pi/2}(A+B|\sin (x+y)|)^{q}C(x)\cos x \,dx\\
&=2mr\int_{0}^{\pi/2}\Big((A+B|\sin (x+y)|)^{q}-(A+B|\sin (x-y)|)^{q}\Big)C(x)\cos x \,dx.
\end{aligned}
$$
For
$$
H(x,y)=\Big((A+B|\sin (x+y)|)^{q}-(A+B|\sin (x-y)|)^{q}\Big)C(x),
$$
we note that
 $$
\begin{cases}
\begin{aligned}
&H(x,y)\leq 0,
\text{ $y\in[0,\frac{\pi}{2}]$~\mbox{and}~$m>1$;}\\
&H(x,y)\geq 0,
\text{ $ y\in[\frac{\pi}{2},\pi]$~\mbox{and}~$m>1$},
\end{aligned}
\end{cases}
$$
and
$$
\begin{cases}
\begin{aligned}
&H(x,y)\geq 0,
\text{ $y\in[0,\frac{\pi}{2}]$~\mbox{and}~$m<1$;}\\
&H(x,y)\leq 0,
\text{ $ y\in[\frac{\pi}{2},\pi]$~\mbox{and}~$m<1$.}
\end{aligned}
\end{cases}
$$
This implies that $\frac{\pi}{2}$ is the maximum of the function $L(y)$ for $m<1$, and $0$ is its maximum for
$m>1.$ The case $m=1$ is trivial and in this case the function $L$ is constant.
\end{proof}

By using the similar arguments as that in the proof of \cite[Lemma 2.2]{David} and Lemma \ref{lem3.2},
we obtain the following:

 \blem\label{lem3.4}
For $m>-1$, $q\geq0$, let
$$
D(r,x)=\int_{-\pi}^{\pi}(A+B|\cos(b-x)|)^{q}(1+r^{2}+2r\cos b)^{m}\,db,
$$
where $A\geq 0$ and $B>0$ are constants.
Then,
$$
D(r,x)\leq
\begin{cases}
\begin{array}{rl}
D(1,0),&
\text{ $m>1;$}\\
D(1,\pi/2),&
\text{ $ m\leq1. $}
\end{array}
\end{cases}
$$
\elem


For a real or complex-valued differentiable function $u=w+iv$ on $\mathbb D$,
the Jacobian matrix  is given by
$$
Du(z)=\begin{pmatrix}
w_x ~~ w_y \\
v_x ~~ v_y \\
\end{pmatrix},
$$
where $z=x+iy$. Its usual operator norm is
$$
|Du(z)|=\sup_{|h|=1} |Du(z)h|=|u_z(z)|+|u_{\overline{z}}(z)|.
$$

The following theorem is a generalization of \cite[Theorem 1.2]{David}.

\bthm\label{thm3.2} Suppose that $u(z)=P_{\alpha,\beta}[f](z)$  on $\mathbb D$ with the complex-valued function $f\in L^{p}(\mathbb{T})$ $(p\geq 1)$. Then, there are a function $B_{\alpha,\beta, p}(r)$  in $(0,1)$ and a constant $B_{\alpha,\beta, p}= \sup_{r\in (0,1)}
B_{\alpha,\beta, p}(r)$ defined in \eqref{eq:495} and \eqref{eq:496} below, so that
$$
|Du(z)|\leq \frac{B_{\alpha,\beta, p}(r)}{(1-r^{2})^{1+\frac{1}{p}}}\|f\|_{ L ^ { p } ( \mathbb{T} )}
\leq \frac{B_{\alpha,\beta, p}}{(1-r^{2})^{1+\frac{1}{p}}}\|f\|_{ L ^ { p } ( \mathbb{T} )}
$$
for $z\in \mathbb{D}$ with $|z|=r$. These inequalities are asymptotically sharp as $\alpha,\beta \rightarrow 0$.
\ethm

\begin{proof}
For a constant $h=e^{i\tau}$,
 we have
$$
|Du(z)h|=\left|\int_{0}^{2\pi} \Big(\big(P_{\alpha,\beta}(z e^{-it})\big)_z e^{i\tau}+\big(P_{\alpha,\beta}(z e^{-it})\big)_{\overline{z}}e^{-i\tau}\Big)
f(e^{it})\,\frac{dt}{2\pi}\right|,
$$
 where $z\in\mathbb D$.
 According to the H\"{o}lder inequality, we obtain
\be\label{eq:6251}
|D u(z)|\leq \|f\|_{p}\max_{\tau}\left(\int_{0}^{2\pi}\Big|\big(P_{\alpha,\beta}(z e^{-it})\big)_z e^{i\tau}+\big(P_{\alpha,\beta}(z e^{-it})\big)_{\overline{z}}e^{-i\tau}\Big|^{q}\frac{dt}{2\pi}\right)^{\frac{1}{q}},
\ee for $q$ such that $\frac{1}{p}+\frac{1}{q}=1$.
Elementary computations show that
$$
\big(P_{\alpha,\beta}(z e^{-it})\big)_z=c_{\alpha,\beta}\frac{(1-|z|^{2})^{\alpha+\beta} \big( e^{-it}(\alpha+1+\beta|z|^{2})-(\alpha+\beta+1)\overline{z}\big)}
{(1-ze^{-it})^{\alpha+2}(1-\bar{z}e^{it})^{\beta+1}}.
$$
Then, for $ z=re^{i\theta}, t=c+\theta,$  
we get
$$
\big(P_{\alpha,\beta}(z e^{-it})\big)_z e^{i\tau}=c_{\alpha,\beta}\frac{(1-r^{2})^{\alpha+\beta}\big(e^{-ic}(\alpha+1+\beta r^{2})-(\alpha+\beta+1)r\big)}
{e^{-i(\tau-\theta)}(1-re^{-ic})^{\alpha+1}(1-re^{ic})^{\beta+2}}.
$$
By making the substitution
 $ e^{ic}=\frac{r+e^{ib}}{1+re^{ib}},$ we have
$$
\big(P_{\alpha,\beta}(z e^{-it})\big)_z e^{-i\tau}=c_{\alpha,\beta}\frac{(-\beta r+(\alpha +1)e^{-ib})(1+re^{-ib})^{\alpha+1}(1+re^{ib})^{\beta+1}}
{e^{-i(\tau-\theta)}(1-r^{2})^{2}}.
$$
Similarly, we have $$
\big(P_{\alpha,\beta}(z e^{-it})\big)_{\overline{z}}e^{-i\tau}=c_{\alpha,\beta}\frac{(-\alpha r+(\beta+1)e^{ib})(1+re^{-ib})^{\alpha+1}(1+re^{ib})^{\beta+1}}
{e^{i(\tau-\theta)}(1-r^{2})^{2}}.
$$


We obtain the following  identity:
$$
\begin{aligned}
&\hspace{-0.5cm}|\big(P_{\alpha,\beta}(z e^{-it})\big)_{\overline{z}}e^{-i\tau}+\big(P_{\alpha,\beta}(z e^{-it})\big)_z e^{i\tau}|^q\\
=&|c_{\alpha,\beta}|^{q}(1-r^{2})^{-2q}(1+r^{2}+2r\cos b)^{\frac{(\alpha+\beta+2)q}{2}}\\
&\left|(-\alpha r+(\beta+1)e^{ib})e^{-i(\tau-\theta)}+(-\beta r+(\alpha +1)e^{-ib})e^{i(\tau-\theta)}\right|^{q}.
\end{aligned}
$$
By combining the above with \eqref{eq:6251}, we obtain
$$
\begin{aligned}
|Du(z)|&\leq\frac{|c_{\alpha,\beta}|\, \|f\|_{p}}{(1-r^{2})^{1+\frac{1}{p}}}\max_{\tau} \left(\int_{0}^{2\pi}(1+r^{2}+2r\cos
b)^{\frac{(\alpha+\beta+2)q}{2}-1}\right. \\
&\times\left. \left|(-\alpha r+(\beta+1)e^{ib})e^{-i(\tau-\theta)}+(-\beta r+(\alpha +1)e^{-ib})e^{i(\tau-\theta)}\right|^{q}
\frac{db}{2\pi}\right)^{\frac{1}{q}}.
\end{aligned}
$$
It is easy to verify that
\beqq
&&\left|(-\alpha r+(\beta+1)e^{ib})e^{-i(\tau-\theta)}+(-\beta r+(\alpha +1)e^{-ib})e^{i(\tau-\theta)}\right|^{q}\\
&=&|(\alpha+\beta+2)\cos(\tau-\theta-b)-(\alpha+\beta)r\cos(\tau-\theta)\\
&&+i(\alpha-\beta)\big(\sin(\tau-\theta-b)+r\sin(\tau-\theta)\big)|^q\\
&\leq&\Big((\alpha+\beta+2)|\cos(\tau-\theta-b)|+|(\alpha+\beta)r\cos(\tau-\theta)|+2|\alpha-\beta|\Big)^q.
\eeqq
Now, for $\rho_1(t)=(t+(\alpha+\beta+2)|\cos(\tau-\theta-b)|+2|\alpha-\beta|)^{q}$, we consider the following inequality:
$$
\rho_1(y)-\rho_1(0)\leq \max_{0\leq t\leq y} |\rho_1{'}(t)|y,
$$
where $y=|(\alpha+\beta)r\cos(\tau-\theta)|.$
Since
$$
\rho_1{'}(t)=q\big(t+(\alpha+\beta+2)|\cos(\tau-\theta-b)|+2|\alpha-\beta|\big)^{q-1},
$$
and
$$
\rho_1(0)=\big((\alpha+\beta+2)|\cos(\tau-\theta-b)|+2|\alpha-\beta|\big)^{q},
$$
it follows that
$$
\begin{aligned}
&\left|(-\alpha r+(\beta+1)e^{ib})e^{-i(\tau-\theta)}+(-\beta r+(\alpha +1)e^{-ib})e^{i(\tau-\theta)}\right|^{q}\\
&\leq \big((\alpha+\beta+2)|\cos(\tau-\theta-b)|+2|\alpha-\beta|\big)^{q}\\
&\hspace{0.5cm}+q\big((\alpha+\beta+2)+|\alpha+\beta|r+2|\alpha-\beta|\big)^{q-1}|\alpha+\beta| r.
\end{aligned}
$$

Let $P(\alpha,\beta,r)=q\big((\alpha+\beta+2)+|\alpha+\beta|r+2|\alpha-\beta|\big)^{q-1}|\alpha+\beta| r$.
We have
$$
\begin{aligned}
|Du(z)|&\leq\frac{|c_{\alpha,\beta}|\|f\|_{p}}{(1-r^{2})^{1+\frac{1}{p}}}
\bigg(P(\alpha,\beta,r)\int_{0}^{2\pi}(1+r^{2}+2r\cos b)^{\frac{(\alpha+\beta+2)q}{2}-1}\frac{db}{2\pi}\\
&\hspace{-0.6cm}+\int_{0}^{2\pi}\big((\alpha+\beta+2)|\cos(\tau-\theta-b)|+2|\alpha-\beta|\big)^{q}(1+r^{2}+2r\cos
b)^{\frac{(\alpha+\beta+2)q}{2}-1}
\frac{db}{2\pi}\bigg)^{\frac{1}{q}}.
\end{aligned}
$$

For $m=(\alpha+\beta+2)q/2-1$ and $\eta=\theta-\tau$, let
$$
L(\eta)=\int_{0}^{2\pi}\big((\alpha+\beta+2)|\cos(\eta-b)|+2|\alpha-\beta|\big)^{q}(1+r^{2}+2r\cos b)^{m}db.
$$
 From Lemma  \ref{lem3.3}, we obtain
$$
V_{\alpha,\beta, p}(r)=\max_{\eta\in[0,2\pi]}L(\eta)=
\begin{cases}
\begin{array}{rl}
L(0),&
\text{ $m>1;$}\\
L(\pi/2),&
\text{ $ m\leq1. $}
\end{array}
\end{cases}
$$

By \eqref{eq:33001},  we may write
\beq\label{eq:495}
&&\hspace{0.5cm}B_{\alpha,\beta, p}(r)=\frac{|c_{\alpha,\beta}|}{(2\pi)^{1-\frac{1}{p}}}
\Big(P(\alpha,\beta,r)F\big(1-\frac{(\alpha+\beta+2)q}{2}, 1-\frac{(\alpha+\beta+2 ) q } {2 } ; 1 ; r ^ { 2 }\big)\\
\nonumber&&\hspace{4.5cm}+V_{p,\alpha,\beta}(r)\Big)^{1/q}
\eeq
and
\be\label{eq:496}
B_{\alpha,\beta, p}=\frac{|c_{\alpha,\beta}|}{(2\pi)^{1-\frac{1}{p}}}
\Big(P(\alpha,\beta,1)\frac{\Gamma\big((\alpha+\beta+2)q-1\big)}{\Gamma^{2}\Big(\frac{(\alpha+\beta+2)q}{2}\Big)}+V_{p,\alpha,\beta}(1)\Big)^{1/q}.
\ee
Then the desired inequalities follow.

For $\alpha=\beta=0$, the function
$B_{\alpha,\beta, p}(r)$ and the constant $B_{\alpha,\beta, p}$ coincide with the corresponding sharp
function and the constant in  \cite[Theorem 1.1]{Kalaj} for  harmonic functions. This implies that our result is asymptotically sharp.
\end{proof}

 We remark that Theorem \ref{thm3.2} improves \cite[Theorem 3.6]{Milos},  in particular, in the case $\alpha=\beta=0$.

\medskip
Next, we estimate the modulus of each first-order partial derivatives of  an ($\alpha,\beta$)-harmonic function
$u$ using the $L^{p}$ norm of the boundary function $f$.

\bthm\label{thm3.3}
Let $u(z) = P_{\alpha,\beta}[f](z)$ be an $(\alpha,\beta)$-harmonic function on $\mathbb D$, where the complex-valued function $f \in L^{p}(\mathbb{T})$
and $p \geq 1$. Then, for $z = re^{i\theta} \in \mathbb{D}$, we have the following:

$(1)$ There exists a function $C_{\alpha,\beta,p}(r)$ such that
$$| u _ { r } (re^{i\theta})| \leq \frac { C _ { \alpha,\beta , p } ( r ) } { ( 1 - r ^ { 2 } ) ^ { 1 + 1 / p } }
\| f \| _ { L ^ { p } ( \mathbb{T} ) } \leq \frac { C _ { \alpha,\beta , p } } { ( 1 - r ^ { 2 } ) ^ { 1 + 1 / p
} } \| f \| _ { L ^ { p } ( \mathbb{T} ) },
$$
where the constant $C_{\alpha,\beta,p}=\sup_{r\in(0,1)}C_{\alpha,\beta,p}(r)$.
The constant $C_{\alpha,\beta,p}$ is asymptotically sharp as $\alpha,\beta\rightarrow 0$.

$(2)$ There exists a function $D_{\alpha,\beta,p}(r)$  such that
$$|u_{\theta}(re^{i\theta})| \leq \frac { D _ { \alpha,\beta,p}(r)}{(1-r^{2})^{1+1/p}}\|f\|_{L^{p}(\mathbb{T})}
\leq\frac{D_{\alpha,\beta,p}}{(1-r^{2})^{1+1/p}}\|f\|_{L^{p} ( \mathbb{T} ) },$$
where the constant $D_{\alpha,\beta,p}=\sup_{r\in(0,1)}
D_{\alpha,\beta,p}(r)$.
If $\alpha=\beta$, then the constant $D_{\alpha,\beta,p}$ is sharp.

$(3)$ There exists a function $E_{\alpha,\beta,p}(r)$  such that
 $$|u_{z}(z)|, | u _ {\overline{z}} (z)| \leq \frac { E _ { \alpha ,\beta, p } ( r ) } { ( 1 - r ^ { 2 } ) ^ { 1
 + 1 / p } } \| f \| _ { L ^ { p } ( \mathbb{T} ) }
 \leq \frac { E_ { \alpha ,\beta, p } } { ( 1 - r ^
 { 2 } ) ^ { 1 + 1 / p } } \|f\| _ { L ^ { p } ( \mathbb{T} ) },
$$
where the constant $E_{\alpha,\beta,p}=
\sup_{r\in(0,1)}E_{\alpha,\beta,p}(r)$.
The constant $E_{\alpha,\beta,p}$ is asymptotically sharp as $\alpha,\beta \rightarrow 0$.

The functions $C_{\alpha,\beta,p}(r), D_{\alpha,\beta,p}(r), E_{\alpha,\beta,p}(r)$ and the constants $C_{\alpha,\beta,p}, D_{\alpha,\beta,p}, E_{\alpha,\beta,p}$ are defined in \eqref{eq:4921}, \eqref{eq:4103}, \eqref{eq:4108}, \eqref{eq:4922}, \eqref{eq:6303}, and \eqref{eq:4109}, respectively.
\ethm

\begin{proof}

 Let $z=re^{i\theta}\in \mathbb D$.

(1) Differentiating both sides of \eqref{eq:4910} with respect to $r$ yields
\be\label{eq:4913}
u_r( r e ^ { i \theta } )=\frac {c _ { \alpha,\beta } } { 2 \pi } ( 1 - r ^ { 2 } ) ^ { \alpha+\beta } \int_{0}^{2\pi }M(r,\theta)f ( e ^ { i t }) dt,
\ee
where  $M(r,\theta)=a(r,\theta)/b(r,\theta)$ with $b(r,\theta)=
{-(1 - r e ^{ i ( \theta - t ) }) ^ {\alpha+2}(1 - r e ^ { i (t- \theta) }) ^ {\beta +2}}$ and
$$a(r,\theta)=2r (\alpha+\beta + 1 )| 1 - r e ^ { i ( \theta - t ) } | ^ { 2 } + ( 1 - r ^ { 2 } )\left( ( \alpha+1
)(e^{i(\theta-t)}-r)+(\beta+1)(e^{i(t-\theta)}-r)\right).
$$

By using the H\"{o}lder inequality,  \eqref{eq:4913} becomes
\be\label{eq:4914}
| u _ { r }( r e ^ { i \theta } ) | \leq |c_{\alpha,\beta}| ( 1 - r ^ { 2 } ) ^ { \alpha+\beta } I _ { 1 } ^ { \frac { 1 } { q
} } \| f \| _ { L ^ { p } ( \mathbb{T} ) },
\ee
where
$$
I_1=\int_{0}^{2\pi}|M(r,\theta)|^{q}\frac{dt}{2\pi}
$$
for $q$ such that $\frac{1}{p}+\frac{1}{q}=1$.
 By using suitable substitutions in \eqref{eq:5281} and \eqref{eq:5282}, we obtain
\be\label{eq:4916}
\begin{aligned}
I_1&
=(1-r^{2})^{1-(\alpha+\beta+2)q}\int_{0}^{2\pi }\frac{|\alpha r+\beta r-(\alpha+1) e^{-is}-(\beta+1) e^{is}|^{q}}  {|1+re^{-is}|^{2-(\alpha+\beta+2)q}}\frac{ds}{2\pi}\\
&=(1-r^{2})^{1-(\alpha+\beta+2)q}\int_{0}^{2\pi }\frac{|\alpha r+\beta r-(\alpha+\beta+2)\cos
s+i(\alpha-\beta)\sin s|^{q}}{|1+re^{-is}|^{2-(\alpha+\beta+2)q}}\frac{ds}{2\pi}.
\end{aligned}
\ee
Let
\be\label{eq:4915}
 \tilde { C} _ { \alpha,\beta , p }(r)=|c_{\alpha,\beta}| \left(\int_{0}^{2\pi }\frac{|\alpha r+\beta
 r-(\alpha+\beta+2)\cos s+i(\alpha-\beta)\sin
 s|^{q}}{|1+re^{-is}|^{2-(\alpha+\beta+2)q}}\frac{ds}{2\pi}\right)^{\frac{1}{q}}.
\ee
It follows from \eqref{eq:4914}, \eqref{eq:4916}, and \eqref{eq:4915} that
\be\label{eq:4919}
| u _ { r }( r e ^ { i \theta } ) | \leq \frac { \tilde { C } _ { \alpha, \beta, p } ( r ) } { ( 1 - r ^ { 2 } ) ^ { 1 + 1 / p }
} \| f \| _ { L ^{p}( \mathbb{T} ) }.
\ee

By using the inequality $\rho_2(x)
\leq \rho_2(0) + x\max_{0\leq \xi\leq x} \rho_2'(\xi)$ with $x>0$ for $\rho_2(t)=(t+|(\alpha+\beta+2)\cos s+i(\alpha-\beta)\sin s|)^{q}$, we have
 $$
\begin{aligned}
&|\alpha r+\beta r-(\alpha+\beta+2)\cos s+i(\alpha-\beta)\sin s|^{q} \\
 \leq& (|\alpha+\beta|r+|(\alpha+\beta+2)\cos s+i(\alpha-\beta)\sin s|)^{q}\\
\leq& |(\alpha+\beta+2)\cos s+i(\alpha-\beta)\sin s|^{q}\\
 &+\max\{q(|\alpha+\beta|r+|(\alpha+\beta+2)\cos s+i(\alpha-\beta)\sin s|)^{q-1}\}|\alpha+\beta|r\\
\leq& \left(|(\alpha+\beta+2)\cos s|+|\alpha-\beta|\right)^{q}\\
&+q(|\alpha+\beta|r+|\alpha+\beta+2|+|\alpha-\beta|)^{q-1}|\alpha+\beta|r.
\end{aligned}
$$
For the integral in the right side of \eqref{eq:4915}, we obtain the estimate
$$
\begin{aligned}
&\int_{0}^{2\pi}\frac{|\alpha r+\beta r-(\alpha+\beta+2)\cos s+i(\alpha-\beta)\sin
s|^{q}}{|1+re^{-is}|^{2-(\alpha+\beta+2)q}}\frac{ds}{2\pi}\\
\leq&\frac{1}{2\pi}\int_{0}^{2\pi}
\frac{\left(|(\alpha+\beta+2)\cos s|+|\alpha-\beta|\right)^{q}}{|1+re^{-is}|^{2-(\alpha+\beta+2)q}}ds\\
&+\frac{q(|\alpha+\beta|r+|\alpha+\beta+2|+|\alpha-\beta|)^{q-1}
|\alpha+\beta|r}{2\pi}\int_{0}^{2\pi}\frac{1}{|1+re^{-is}|^{2-(\alpha+\beta+2)q}}ds\\
:=&\frac{1}{2\pi}I_{11}+\frac{q(|\alpha+\beta|r+|\alpha+\beta+2|+|\alpha-\beta|)^{q-1}|\alpha+\beta|r}{2\pi}I_{12}.
\end{aligned}
$$
We rewrite $I_{11}$ as
$$
\begin{aligned}
G(r,0):=I_{11}&=\int_{0}^{2\pi}\frac{\left(|(\alpha+\beta+2)\cos s|+|\alpha-\beta|\right)^{q}}{|1+re^{-is}|^{2-(\alpha+\beta+2)q}}ds\\
 &\leq\int_{0}^{2\pi}\left(|(\alpha+\beta+2)\cos s|+|\alpha-\beta|\right)^{q}(1+r^{ 2 } + 2 r \cos s )^{\frac { (\alpha+\beta+2)q}{2}-1}d s.
\end{aligned}
$$
For $(\alpha+\beta+2)q\leq 4$, that is $\frac { ( \alpha+\beta+2 ) q  } { 2 }-1<1$, according to Lemma \ref{lem3.4}, we have
$$
\begin{aligned}
G ( r ,0) &\leq G (1, \frac{\pi}{2}) \\
&= \int_{0}^{2\pi}\left(|(\alpha+\beta+2)\cos (s-\frac{\pi}{2})|+|\alpha-\beta|\right)^{q}(2+2\cos s)^{\frac{(\alpha +\beta+2)q}{2}-1}ds\\
&= 2^{\frac{(\alpha+\beta+2) q}{2}-1}\int _ { 0 }^{2\pi}\big(|(\alpha+\beta+2)\sin s|+|\alpha-\beta|\big)^{q} ( 1+ \cos s )^{ \frac{(\alpha +\beta+2) q}{2}-1}ds.
\end{aligned}
$$
If $(\alpha+\beta+2)q> 4$, then $\frac { (\alpha+\beta+2 ) q  } { 2 }-1> 1$. From Lemma \ref{lem3.4}, it again follows that
$$
G ( r ,0) \leq G \big( 1 ,0\big) = \int _ { 0 } ^ { 2 \pi }\big(|(\alpha+\beta+2)\cos s|+|\alpha-\beta|\big)^{q} (2+2 \cos s ) ^{ \frac{(\alpha +\beta+2) q}{2}-1} d s.
$$

Let
\be\label{eq:4921}
C_{\alpha,\beta,p}(r)=|c_{\alpha,\beta}|\left(\frac{1}{2\pi}I_{11}
+\frac{q(|\alpha+\beta|r+|\alpha+\beta+2|+|\alpha-\beta|)^{q-1}|\alpha+\beta|r}{2\pi}I_{12}\right)^{\frac{1}{q}}.
\ee
Then, by \eqref{eq:33001},
\be\label{eq:4922}
C_{\alpha,\beta,p}=
\begin{cases}
\begin{array}{rl}
|c_{\alpha,\beta}|\Big(c_1(1)\frac{\Gamma\big((\alpha+\beta+2)q-1\big)}
{\Gamma^2\big(\frac{(\alpha+\beta+2)q}{2}\big)}
+\frac{1}{2\pi}G\big(1,\frac{\pi}{2}\big)\Big)^{\frac{1}{q}},&
\text{ $(\alpha+\beta+2)q\leq4;$}\\
|c_{\alpha,\beta}|\Big(c_1(1)\frac{\Gamma\big((\alpha+\beta+2)q-1\big)}{\Gamma^2
\big(\frac{(\alpha+\beta+2)q}{2}\big)}+\frac{1}{2\pi}G\big(1,0\big)\Big)^{\frac{1}{q}},&
\text{ $(\alpha+\beta+2)q> 4,$}\\
\end{array}
\end{cases}
\ee
where
$c_1(r)=q(|\alpha+\beta|r+|\alpha+\beta+2|+|\alpha-\beta|)^{q-1}|\alpha+\beta|r$.
Then $\tilde { C } _ { \alpha,\beta , p } ( r ) \leq C_{\alpha, \beta, p}(r) \leq C_{\alpha, \beta, p}$.
Using  \eqref{eq:4919} yields the required inequalities.

Specially, if $\alpha=\beta=0$, then
$$
C_{0,0,p}(r)=2\left(\frac{1}{2\pi}\int_{0}^{2\pi}|\cos s|^{q}|1+re^{-is}|^{2q-2}ds\right)^{\frac{1}{q}}
$$
and
$$
C_{0,0,p}=\frac{4^{\frac{1}{p}}}{\pi^{\frac{1}{q}}}\left(\int_{0}^{2\pi}|\cos s|^{q}(1+\cos s)^{q-1}ds\right)^{\frac{1}{q}}.
$$
By using the similar arguments as that in the proof of \cite[Theorem 1.2]{Yong}, we see that $C_{0,0,p}$ is sharp.

\smallskip

(2) Differentiating both sides of the formula \eqref{eq:4910} with respect to $\theta$, we have
\be\label{eq:4101}
u_{\theta}(re^{i\theta})
=-\frac{c_{\alpha,\beta}}{2\pi}(1-r^{2})^{\alpha+\beta+1}r\int_{0}^{2\pi}N(r,\theta)f(e^{it})dt,
\ee
where
$$
N(r,\theta)=\frac{i[(\alpha+1)(r-e^{i(\theta-t)})+(\beta+1)(e^{i(t-\theta)}-r)]}
{(1-re^{i(\theta-t)})^{\alpha+2}(1-re^{i(t-\theta)})^{\beta+2}}.
$$
 Then applying the H\"{o}lder inequality yields
$$
|u_{\theta}(re^{i\theta})|\leq |c_{\alpha,\beta}|r(1-r^{2})^{\alpha+\beta+1}I_2^{\frac{1}{q}}\|f\|_{L^{p}(\mathbb{T})},
$$
where
$$
I_2=\int_{0}^{2\pi}|N(r,\theta)| ^{q}\frac{dt}{2\pi}
$$
and $q$ is such that $\frac{1}{p}+\frac{1}{q}=1$.

By using  substitutions \eqref{eq:5281} and \eqref{eq:5282}, we have
$$
\begin{aligned}
I_2&=\int_{0}^{2\pi}\frac{|(\alpha+1)(r-\frac{1+re^{is}}{r+e^{is}})+(\beta+1)(\frac{r+e^{is}}{1+re^{is}}-r)|^{q}}
{\big(\frac{1-r^{2}}{|1+re^{-is}|}\big)^{(\alpha+\beta+4)q}} \frac{1-r^{2}}{|1+re^{is}|^{2}}\frac{ds}{2\pi}\\
&=(1-r^{2})^{1-(\alpha+\beta+3)q}\int_{0}^{2\pi}\frac{|e^{is}(\beta+1)-(\alpha+1)e^{-is}+(\beta-\alpha)r|^{q}}
{|1+re^{-is}|^{2-(\alpha+\beta+2)q}}\frac{ds}{2\pi}.
\end{aligned}
$$
By the inequality $\rho_3(x)
\leq \rho_3(0) + x\max_{0\leq\xi\leq x} \rho_3'(\xi)$ with $x>0$ for $\rho_3(t)=\big(|(\beta-\alpha)\cos s+i(\alpha+\beta+2)\sin s|+t\big)^{q}$, we obtain
$$
\begin{aligned}
&|e^{is}(\beta+1)-(\alpha+1)e^{-is}+(\beta-\alpha)r|^{q}\\
\leq& \big(|(\beta-\alpha)\cos s+i(\alpha+\beta+2)\sin s|+|\beta-\alpha|r\big)^{q}\\
\leq& |(\beta-\alpha)\cos s+i(\alpha+\beta+2)\sin s|^{q}\\
&+q \max \{\left(|(\beta-\alpha)\cos s+i(\alpha+\beta+2)\sin s|+|\beta-\alpha|r\right)^{q-1}\}|\beta-\alpha|r\\
\leq &(|\beta-\alpha|+|(\alpha+\beta+2)\sin s|)^{q}\\
&+q\left(|\beta-\alpha|+|\alpha+\beta+2|+|\beta-\alpha|r\right)^{q-1}|\beta-\alpha|r.
\end{aligned}
$$
We have the following estimate:
$$
\begin{aligned}
\int_{0}^{2\pi}&\frac{|e^{is}(\beta+1)-(\alpha+1)e^{-is}+(\beta-\alpha)r|^{q}}{|1+re^{-is}|^{2-(\alpha+\beta+2)q}}\frac{ds}{2\pi}\\
\leq& \frac{1}{2\pi}I_{13}
+ q\frac{\left(|\beta-\alpha|+|\alpha+\beta+2|+|\beta-\alpha|r\right)^{q-1}|\beta-\alpha|r}{2\pi}I_{12},
\end{aligned}
$$
where
$$
\begin{aligned}
I_{13}=&\int_{0}^{2\pi}\frac{(|\beta-\alpha|+|(\alpha+\beta+2)\sin s|)^{q}}{|1+re^{-is}|^{2-(\alpha+\beta+2)q}}\frac{ds}{2\pi}\\
=&\int_{0}^{2\pi}(|\beta-\alpha|+|(\alpha+\beta+2)\sin s|)^{q}(1+r^{ 2 } + 2 r \cos s )^{\frac { ( \alpha+\beta+2 ) q}{2}-1}d s:=G(r,\frac{\pi}{2}).
\end{aligned}
$$

Let
\be\label{eq:4103}
\begin{aligned}
D_{\alpha,\beta,p}(r)&=|c_{\alpha,\beta}|r\left(\frac{1}{2\pi}I_{13}+
q\frac{\left(|\beta-\alpha|+|\alpha+\beta+2|+|\beta-\alpha|r\right)^{q-1}|\beta-\alpha|r}{2\pi}I_{12}\right)^{\frac{1}{q}}.
\end{aligned}
\ee
By Lemma \ref{lem3.4}, we have
\be\label{eq:6303}
D_{\alpha,\beta,p}=
\begin{cases}
\begin{array}{rl}
|c_{\alpha,\beta}|\left(\frac{1}{2\pi}G\big(1,\frac{\pi}{2})+c_2(r)\frac{\Gamma\big((\alpha+\beta+2)q-1\big)}{\Gamma^{2}
\big(\frac{(\alpha+\beta+2)q}{2}\big)}\right)^{\frac{1}{q}},&
\text{$(\alpha+\beta+2)q\leq4$};\\
|c_{\alpha,\beta}|\left(\frac{1}{2\pi}G(1,0\big)
+c_2(r)\frac{\Gamma\big((\alpha+\beta+2)q-1\big)}{\Gamma^{2}
\big(\frac{(\alpha+\beta+2)q}{2}\big)}\right)^{\frac{1}{q}},&
\text{$(\alpha+\beta+2)q> 4,$}\\
\end{array}
\end{cases}
\ee
where
$c_2(r)=q\left(|\beta-\alpha|+|\alpha+\beta+2|+|\beta-\alpha|r\right)^{q-1}|\beta-\alpha|r$.
Thus, it holds that
$$|u_{\theta}( r e ^ { i \theta } )|\leq\frac{D_{\alpha,\beta,p}(r)}{(1-r^{2})^{1+\frac{1}{p}}}\|f\|_{L^{p}(\mathbb{T})}\leq \frac{D_{\alpha,\beta,p}}{(1-r^{2})^{1+\frac{1}{p}}}\|f\|_{L^{p}(\mathbb{T})}.$$

For $\alpha=\beta$,  by using \eqref{eq:6301} of Lemma \ref{lem3.4}, we have
$$
\begin{aligned}
D_{\alpha,\alpha,p}(r)=&\frac{\Gamma^{2}(\alpha+1)}{|\Gamma(2\alpha+1)|}
\frac{(2\alpha+2)r}{\pi^{\frac{1}{q}}}\left(\int_{0}^{\pi}(\sin s)^{q}(1+r^{2}-2r\cos s)^{(\alpha+1)q-1}ds\right)^{\frac{1}{q}}\\
=&\frac{\Gamma^{2}(\alpha+1)}{|\Gamma(2\alpha+1)|}\frac{(2\alpha+2)r}{\pi^{\frac{1}{q}}}\\
&\times\left(B\Big(\frac{1+q}{2},\frac{1}{2}\Big)F\Big(1-(\alpha+1)q,1-\big(\alpha+\frac{3}{2}\big)q;1+\frac{q}{2};r^{2}\Big)\right)^{\frac{1}{q}}
\end{aligned}
$$
and
\beq
\nonumber \hspace{-1cm} D_{\alpha,\alpha,p}&=&\frac{\Gamma^{2}(\alpha+1)}{|\Gamma(2\alpha+1)|} \frac{2\alpha+2}{\pi^{\frac{1}{q}}}\\
&&\nonumber \times \left(B\Big(\frac{1+q}{2},\frac{1}{2}\Big)F\Big(1-(\alpha+1)q,1-\big(\alpha+\frac{3}{2}\big)q;1+\frac{q}{2};1\Big)\right)^{\frac{1}{q}}.
\eeq
 From \cite[Theorem 1.2]{Yong}, we see that $D_{\alpha,\alpha,p}$ is sharp.

\smallskip
(3) Differentiating on both sides of the formula \eqref{eq:4910} with respect to $z$ yields
\be\label{eq:41016}
\begin{aligned}
\nonumber u_z( z ) &= c _ { \alpha,\beta } ( 1 - r ^ { 2 } ) ^ { \alpha+\beta } \\
\nonumber&\int_ { 0 } ^ { 2 \pi } \frac { - ( \alpha+\beta + 1 ) \overline { z } ( 1 - re ^ {i (\theta-t) } ) + ( 1 +
\alpha) ( 1 - r ^ { 2 } ) e ^ { - i t }  } { ( 1 - re ^ {i (\theta-t) } ) ^ { \alpha + 2 }( 1 - re ^ {i
(t-\theta) } ) ^ { \beta + 1 } } f ( e ^ { i t } ) \frac { d t } { 2 \pi }.
\end{aligned}
\ee
Then, by the H\"{o}lder inequality, we have
\be\label{eq:4107}
|u_{z}(z)| \leq |c _ {\alpha,\beta}|( 1 - r ^ { 2 } ) ^ { \alpha+\beta } I _ { 3 } ^ { \frac { 1 } {
q} } | | f | | _ { L ^ {p} (\mathbb{ T }) },
\ee
where
$$I_3=\int _ { 0 } ^ { 2 \pi } \frac {| - ( \alpha+\beta + 1 ) \overline { z } ( 1 - re ^ {i (\theta-t) } )
+ ( 1 + \alpha) ( 1 - r ^ { 2 } ) e ^ { - i t }|^{q}  } { | 1 - re ^ {i (\theta-t) } | ^ { (\alpha +
\beta+3)q } }\frac { d t } { 2 \pi },$$
where $q$ is such that $ \frac { 1 } { p } + \frac { 1 } { q } = 1 $.

After the change of variables  \eqref{eq:5281} and \eqref{eq:5282}, we have
$$
\begin{aligned}
I_3&=\int _ { 0 } ^ { 2 \pi } \frac {| - ( \alpha+\beta + 1 )
(r\frac{r+e^{is}}{1+re^{is}})\frac{e^{is}(1-r^{2})}{r+e^{is}} + ( 1 + \alpha) ( 1 - r ^ { 2 } )|^{q}  } {(\frac{1-r^{2}}{| 1+re ^ {-is}|}) ^ { (\alpha + \beta+3)q } }\frac{1-r^{2}}{| 1+re ^
{is}|^{2}}\frac { d s } { 2 \pi }\\
&=(1-r^{2})^{1-(\alpha+\beta+2)q} \int _ { 0 } ^ { 2 \pi }\frac{|(\alpha+1)-\beta re^{is}|^{q}}
{|1+re ^ {-is}|^{2-(\alpha+\beta+2)q}}\frac { d s } { 2 \pi }.
\end{aligned}
$$
Let
$$
\tilde { E }_{\alpha,\beta,p}(r)=|c _ { \alpha,\beta }|\left(\int _ { 0 } ^ { 2 \pi }\frac{|(\alpha+1)-\beta
re^{is}|^{q}}{|1+re ^ {-is}|^{2-(\alpha+\beta+2)q}}\frac { d s } { 2 \pi }\right)^{\frac{1}{q}}.
$$
Then \eqref{eq:4107} becomes
\be\label{eq:41011}
| u _ { z } (z)| \leq \frac { \tilde {E} _ { \alpha,\beta , p } ( r ) } { ( 1 - r ^ { 2 } ) ^ { 1 + 1 / p }
} | | f | _ { L ^ {p} (\mathbb{ T} ) }.
\ee
Since $|(\alpha+1)-\beta re^{is}|^{q}\leq (|\alpha+1|+|\beta|r)^{q}$, it follows that
\be\label{eq:4108}
\tilde {E}_{\alpha,\beta,p}(r)\leq |c_{\alpha,\beta}|(|\alpha+1|+|\beta|r)\Big(\frac{1}{2\pi}\int _ { 0 }
^ { 2 \pi }\frac{ds}{|1+re ^ {-is}|^{2-(\alpha+\beta+2)q}} \Big)^{\frac{1}{q}}:=E_{\alpha,\beta,p}(r).
\ee
By \eqref{eq:33001}, we can rewrite $E_{\alpha,\beta,p}(r)$ as
$$
E_{\alpha,\beta,p}(r)=|c_{\alpha,\beta}|
(|\alpha+1|+|\beta|r)\left(F\Big(1-\frac{(\alpha+\beta+2)q}{2},1-\frac{(\alpha+\beta+2)q}{2};1;r^{2}\Big)\right)^{\frac{1}{q}}.\\
$$
and
\be\label{eq:4109}
E_{\alpha,\beta,p}=|c_{\alpha,\beta}|(|\alpha+1|+|\beta|)
\left(\frac{\Gamma\big((\alpha+\beta+2)q-1\big)}{\Gamma^{2}\Big(\frac{(\alpha+\beta+2)q}{2}\Big)}\right)^{\frac{1}{q}}.
\ee
Then, it follows from \eqref{eq:41011} and \eqref{eq:4108} that the desired inequalities hold with $E_{\alpha,\beta,p}(r)\leq E_{\alpha,\beta,p}$.

 If $\alpha=\beta=0$, then
 $$E_{0,0,p}(r)=|c_{\alpha,\beta}|
(|\alpha+1|+|\beta|r)\left(F\Big(1-q,1-q;1;r^{2}\Big)\right)^{\frac{1}{q}}$$ and
$$E_{0,0,p}=\big( \frac { \Gamma ( 2 q - 1 ) } {  \Gamma^ { 2 } ( q )  } \big) ^ { \frac {
1 } { q } }.$$ 
It can be confirmed that $E_{0,0,p}$ is sharp, the details  are in the proof of \cite[Theorem 1.2]{Yong}.

The estimate for $|u_{\overline{z}}(z)|$ is similar.
We omit the proof.
\end{proof}

\medskip

By using the $L^{p}$ norm of the boundary function $f$, we estimate the $ M_{p}(r,\cdot)$ of each first-order partial derivative of
the ($\alpha,\beta$)-harmonic function $u(z)=P_{\alpha,\beta}[f](z)$.

\bthm\label{thm3.4}
Let $u(z) = P_{\alpha,\beta}[f](z)$ be an $(\alpha,\beta)$-harmonic function on $\mathbb D$ with the complex-valued function $f \in L^{p}(\mathbb{T})$,
where $p \geq 1$. Then, for $z = re^{i\theta} \in \mathbb{D}$, the following
statements hold:

$(1)$ There exists a function $A_{\alpha,\beta}(r)$
such that
$$M _ { p } ( r , u _ { r } ) \leq \frac { A_ { \alpha,\beta } ( r ) } { 1 - r ^ { 2 } } \| f \| _ { L ^
{ p } (\mathbb T) } \leq \frac {A_ { \alpha,\beta } } { 1 - r ^ { 2 } } \| f \|  _ { L ^ { p } ( \mathbb{T} )},
$$
where $A_{\alpha,\beta}=\sup_{r\in(0,1)}A_{\alpha,\beta}(r)$.
The constant $A_{\alpha,\beta}$ is asymptotically sharp as $\alpha,\beta\rightarrow 0$.

$(2)$ There exists a function $B_{\alpha,\beta}(r)$ such that
$$
M _ { p } ( r , u _ {\theta} ) \leq \frac {B_ { \alpha,\beta } ( r ) } { 1 - r ^ { 2 } }
 \|f\|_{ L^{p}(\mathbb T ) } \leq \frac{B_ { \alpha,\beta } }{ 1 - r ^ { 2 } } \| f \|  _ { L ^ { p } ( \mathbb{T} )},
$$
where $B_{\alpha,\beta}=\sup_{r\in(0,1)}
B_{\alpha,\beta}(r)$.
The constant $B_{\alpha,\beta}$ is asymptotically sharp as $\alpha,\beta\rightarrow 0$.

$(3)$ There exists a function $C_{\alpha,\beta}(r)$  such that
$$M _ { p } ( r , u _ {z} ) \leq \frac {C_ { \alpha,\beta } ( r ) } { 1 - r ^ { 2 } } \| f \| _ { L ^ { p
} (\mathbb T ) } \leq \frac {C_ { \alpha,\beta } } { 1 - r ^ { 2 } } \|f\|_{L^{p}(\mathbb{T})},
$$
where $C_{\alpha,\beta}=\sup_{r\in(0,1)}
C_{\alpha,\beta}(r)$.
The constant $C_{\alpha,\beta}$ is asymptotically sharp as $\alpha,\beta \rightarrow 0$.

The functions $A_{\alpha,\beta}(r), B_{\alpha,\beta}(r)$, $C_{\alpha,\beta}(r)$, and the constants $A_{\alpha,\beta}, B_{\alpha,\beta}$, $C_{\alpha,\beta}$ are defined in \eqref{eq:41022}, \eqref{eq:6111}, \eqref{eq:41020}, \eqref{eq:41023}, \eqref{eq:6302}, and \eqref{eq:41021}, respectively.\ethm

\begin{proof}

 Let $z=re^{i\theta}$ be in $ \mathbb D$.

(1) According to equation \eqref{eq:4913}, we have
$$|u_{r}(re^{i\theta})|\leq |c_{\alpha,\beta}|(1-r^{2})^{\alpha+\beta}\int_{0}^{2\pi}|M(r,\theta)
f(e^{it})|\frac{dt}{2\pi},$$
where $M(r,\theta)$ is defined as in the proof of Theorem \ref{thm3.3}.
Then, by Jensen's inequality, we have
$$
\begin{aligned}
|u_{r}(re^{i\theta})|^{p}&\leq\left(|c_{\alpha,\beta}|(1-r^{2})^{\alpha+\beta}\right)^{p}I_4^{p}\left(\int_{0}^{2\pi}
\frac{|M(r,\theta)|}{I_4}|f(e^{it})|\frac{dt}{2\pi}\right)^{p}\\
&\leq (|c_{\alpha,\beta}|(1-r^{2})^{\alpha+\beta})^{p}I_4^{p-1}\int_{0}^{2\pi}|M(r,\theta)|
|f(e^{it})|^{p}\frac{dt}{2\pi},
\end{aligned}
$$
where
$$
I_4=\int_{0}^{2\pi}|M(r,\theta)|
\frac{dt}{2\pi}.
$$
Integrating both sides of the above inequality yields
$$
\begin{aligned}
&\frac{1}{2\pi}\int_{0}^{2\pi}|u_{r}(re^{i\theta})|^{p}d\theta\\
& \leq \big(|c_{\alpha,\beta}|(1-r^{2})^{\alpha+\beta}\big)^{p}I_4^{p-1} \int_{0}^{2\pi}\int_{0}^{2\pi}|M(r,\theta)|
|f(e^{it})|^{p}\frac{dt}{2\pi}\frac{d\theta}{2\pi}\\
&=\big(|c_{\alpha,\beta}|(1-r^{2})^{\alpha+\beta}\big)^{p}I_4^{p}\|f\|_{L^{p}(\mathbb{T})} ^ { p }.
\end{aligned}
$$
By change of variables in \eqref{eq:5281} and \eqref{eq:5282}, we see that
$$
I_4=(1-r^{2})^{-1-\alpha-\beta}\int_{0}^{2\pi}\frac{|\alpha r+\beta r-(\alpha+\beta+2)\cos s+i(\alpha-\beta)\sin s|}{|1+re^{-is}|^{-(\alpha+\beta)}}\frac{ds}{2\pi}.
$$
Hence,
\be \label{eq:41014}
M_{p}(r,u_{r})\leq\frac{\tilde{A}_{ \alpha,\beta } ( r ) }{1-r^{2}} \| f \| _ { L ^ { p } ( \mathbb{T} ) },
\ee
where
$$
\tilde { A} _ { \alpha,\beta } ( r ) =|c_{\alpha,\beta}|\int_{0}^{2\pi}\frac{|\alpha r+\beta r-(\alpha+\beta+2)\cos s+i(\alpha-\beta)\sin s|}{|1+re^{-is}|^{-(\alpha+\beta)}}\frac{ds}{2\pi}.
$$
Let
\be\label{eq:41022}
\begin{aligned}
A_{\alpha,\beta}(r)&=|c_{\alpha,\beta}|\left(\frac{|\alpha+\beta|r}{2\pi}\int_{0}^{2\pi}|1+re^{-is}|^{\alpha+\beta}ds
+\frac{\alpha+\beta+2}{2\pi}\right.\\
&\left. \hspace{0.5cm}\times\int_{0}^{2\pi}|\cos s||1+re^{-is}|^{\alpha+\beta}ds
+\frac{|\alpha-\beta|}{2\pi}\int_{0}^{2\pi}|\sin s||1+re^{-is}|^{\alpha+\beta}ds\right)
\end{aligned}
\ee
and
\be\label{eq:41023}
A_{\alpha,\beta}=
\begin{cases}
\begin{aligned}
&|c_{\alpha,\beta}|\left(\frac{\alpha+\beta+2+|\alpha-\beta|}{\pi}2^{\frac{\alpha+\beta}{2}-1}
\int_{0}^{2\pi}|\cos s|(1+\cos s)^{\frac{\alpha+\beta}{2}}ds\right.\\
&\hspace{1cm}\left.+|\alpha+\beta|\frac{\Gamma(\alpha+\beta+1)}{\Gamma^{2}\big(\frac{\alpha+\beta+2}{2}\big)}\right),
\quad\text{ $\alpha+\beta\geq 2;$}\\
&|c_{\alpha,\beta}|\left(\frac{\alpha+\beta+2+|\alpha-\beta|}{\pi}2^{\frac{\alpha+\beta}{2}-1}\int_{0}^{2\pi}|\sin s|(1+\cos s)^{\frac{\alpha+\beta}{2}}ds\right.\\
&\hspace{1cm}\left.+|\alpha+\beta|\frac{\Gamma(\alpha+\beta+1)}{\Gamma^{2}\big(\frac{\alpha+\beta+2}{2}\big)}\right),
\quad\text{ $-1<\alpha+\beta<2. $}
\end{aligned}
\end{cases}
\ee

By using Lemma \ref{lem3.4}, we have that $ \tilde{A}_{\alpha,\beta}(r) \leq A_{\alpha,\beta}(r) \leq A _ { \alpha,\beta}$. Therefore, by \eqref{eq:41014}, the desired inequalities  follow.

Specially, if $\alpha=\beta=0$, then $A_{0,0}(r)=A_{0,0}=\frac{1}{\pi}\int_{0}^{2\pi}|\cos s|ds=\frac{4}{\pi}$.
The process of verifying that $A_{0,0}$ is sharp can be found in the proof of \cite[Theorem 1.3]{Yong}.

\smallskip

(2) Equation \eqref{eq:4101} gives that
$$
\begin{aligned}
|u_{\theta}(re^{i\theta})|&\leq \frac{|c_{\alpha,\beta}|}{2\pi}(1-r^{2})^{\alpha+\beta+1}\\
&\hspace{0.5cm}\times\int_{0}^{2\pi}\frac{r|(\alpha+1)(r-e^{i(\theta-t)})+(\beta+1)(e^{i(t-\theta)}-r)|}
{(1+r^{2}-2r\cos(\theta-t))^{2+\frac{\alpha+\beta}{2}}}|f(e^{it})|\,dt.
\end{aligned}
$$
By using the similar arguments as that in (1) for $|u_{r}(re^{i\theta})|$, we get
$$
\begin{aligned}
|u_{\theta}(re^{i\theta})|^{p}\leq& \left(|c_{\alpha,\beta}|(1-r^{2})^{\alpha+\beta+1}\right)^{p}I_5^{p-1}\\
&\times\int_{0}^{2\pi}\frac{r|(\alpha+1)(r-e^{i(\theta-t)})+(\beta+1)(e^{i(t-\theta)}-r)|}
{(1+r^{2}-2r\cos(\theta-t))^{2+\frac{\alpha+\beta}{2}}}|f(e^{it})|^{p}\frac{dt}{2\pi},
\end{aligned}
$$
where
$$
I_5=\int_{0}^{2\pi}\frac{r|(\alpha+1)(r-e^{i(\theta-t)})+(\beta+1)(e^{i(t-\theta)}-r)|}
{(1+r^{2}-2r\cos(\theta-t))^{2+\frac{\alpha+\beta}{2}}}\,dt.
$$
Consequently,
\be\label{eq:5231}
\frac{1}{2\pi}\int_{0}^{2\pi}|u_{\theta}(re^{i\theta})|^{p}d\theta \leq (|c_{\alpha,\beta}|(1-r^{2})^{\alpha+\beta+1})^{p}I_5^{p}\|f\|_{L^{p}(\mathbb{T})} ^ { p }.
\ee
By using \eqref{eq:5281} and \eqref{eq:5282} again, we have
\be\label{eq:5232}
I_5=(1-r^{2})^{-2-\alpha-\beta}\int_{0}^{2\pi}\frac{r|e^{is}(\beta+1)-(\alpha+1)e^{-is}+(\beta-\alpha)r|}{|1+re^{-is}|^{-\alpha-\beta}}
\frac{ds}{2\pi}.
\ee
Let
\be\label{eq:6111}
\begin{aligned}
B_{\alpha,\beta}(r)&=|c_{\alpha,\beta}|\left(\frac{|\beta-\alpha|r}{2\pi}\int_{0}^{2\pi}|1+re^{-is}|^{\alpha+\beta}ds
+\frac{\alpha+\beta+2}{2\pi}\right.\\
&\left.\int_{0}^{2\pi}|\sin s||1+re^{-is}|^{\alpha+\beta}ds
+\frac{|\beta-\alpha|}{2\pi}\int_{0}^{2\pi}|\cos s||1+re^{-is}|^{\alpha+\beta}ds\right).
\end{aligned}
\ee

Therefore, it follows from \eqref{eq:5231} and \eqref{eq:5232} that
$$
M_{p}(r,u_{\theta})\leq\frac{B_{\alpha,\beta}(r)}{1-r^{2}}\|f\|_{L^{p}(\mathbb{T})}.
$$
By using Lemma \ref{lem3.4}, we can directly verify that
\be\label{eq:6302}
B_{\alpha,\beta}=
\begin{cases}
\begin{aligned}
&|c_{\alpha,\beta}|\left(\frac{\alpha+\beta+2+|\beta-\alpha|}{\pi}2^{\frac{\alpha+\beta}{2}-1}\int_{0}^{2\pi}|\cos s|(1+\cos s)^{\frac{\alpha+\beta}{2}}ds\right. \\
&\hspace{1cm}\left.+|\beta-\alpha|\frac{\Gamma(\alpha+\beta+1)}{\Gamma^{2}(\frac{\alpha+\beta+2}{2})}\right),
\quad\text{ $\alpha+\beta\geq 2;$}\\
&|c_{\alpha,\beta}|\left(\frac{\alpha+\beta+2+|\beta-\alpha|}{\pi}2^{\frac{\alpha+\beta}{2}-1}\int_{0}^{2\pi}|\sin s|(1+\cos s)^{\frac{\alpha+\beta}{2}}ds\right.\\
&\hspace{1cm}\left.+|\beta-\alpha|\frac{\Gamma(\alpha+\beta+1)}{\Gamma^{2}(\frac{\alpha+\beta+2}{2})}\right),
\quad \text{ $-1<\alpha+\beta<2.$}
\end{aligned}
\end{cases}
\ee

 When $\alpha=\beta=0$, we have
$$
B_{0,0}(r)=\int_{0}^{2\pi}r|e^{is}-e^{-is}|\frac{ds}{2\pi}=4\int_{0}^{\pi}r\sin s \frac{ds}{2\pi}.
$$
Obviously,
$$
B_{0,0}(r)\leq\frac{4}{\pi}=B_{0,0}.
$$
As in the proof of \cite[Theorem 1.3]{Yong}, it can be verified that $B_{0,0}$ is sharp.

(3) It follows from \eqref{eq:41016} that
\be\label{eq:41017}
\begin{aligned}
|u_{z}(z)| &\leq |c_{\alpha,\beta}| ( 1 - r ^ { 2 } ) ^ {\alpha+\beta}\\
& \int_{0}^{2\pi} \frac {|- (\alpha+\beta + 1 ) \overline { z } ( 1 - re ^ {i (\theta-t) } ) + ( 1 + \alpha) ( 1 - r ^ { 2 } ) e ^ { -
i t } | } { |1 - re ^ {i (\theta-t) } | ^ { \alpha + \beta + 3 } } |f ( e ^ { i t } )| \frac { d t } { 2\pi }.
\end{aligned}
\ee
By Jensen's inequality, \eqref{eq:41017} leads to
$$
\begin{aligned}
|u_{z}(z)|^{p} \leq& \big(|c_{\alpha,\beta}|( 1 - r ^ { 2 } ) ^ {\alpha+\beta}\big)^{p}\\
&\times\left(I_6\cdot \int _ { 0 } ^ {2 \pi}\frac {|- ( \alpha+\beta + 1 ) \overline { z } ( 1 - re ^ {i (\theta-t) } ) + ( 1 + \alpha) ( 1 - r ^ { 2
} ) e ^ { - i t } | } { |1 - re ^ {i (\theta-t) } | ^ { \alpha + \beta + 3 }I_6 }|f ( e ^ { i t } )| \frac
{ d t } { 2 \pi }\right)^{p}\\
\leq& \big(|c_{\alpha,\beta}|( 1 - r ^ { 2 } ) ^ {\alpha+\beta}\big)^{p}I_6^{p-1}\\
&\times\int _ { 0 } ^ {2 \pi} \frac {|- (\alpha+\beta + 1 ) \overline {z}(1-re^{i(\theta-t)})+(1+\alpha)(1-r^{2})e^{ -
i t } | } { |1 - re ^ {i (\theta-t) } | ^ { \alpha + \beta + 3 } }|f ( e ^ { i t } )|^{p} \frac { d t } {
2 \pi },
\end{aligned}
$$
where
$$
I_6= \int _ { 0 } ^ {2 \pi} \frac {|- ( \alpha+\beta + 1 ) \overline { z } ( 1 - re ^ {i (\theta-t) } ) +
( 1 + \alpha) ( 1 - r ^ { 2 } ) e ^ { - i t } | } { |1 - re ^ {i (\theta-t) } | ^ { \alpha + \beta + 3 }
}\frac { d t } { 2 \pi }.
$$
After the change of variables as in \eqref{eq:5281} and \eqref{eq:5282}, we have
$$
I_6=(1-r^{2})^{-1-(\alpha+\beta)} \int _ { 0 } ^ { 2 \pi }\frac{|(\alpha+1)-\beta re^{is}|}
{|1+re ^ {-is}|^{-(\alpha+\beta)}}\frac { d s } { 2 \pi }.
$$
It follows that
\be\label{eq:41018}
\begin{aligned}
&\frac{1}{2\pi}\int _ { 0 } ^ {2 \pi}|u_{z}(z)|^{p}d\theta \leq (|c_{\alpha,\beta}|( 1 - r ^ { 2 } ) ^ {\alpha+\beta})^{p}I_6^{p-1}\cdot \\
& \int _ { 0 } ^ {2 \pi}\left(\int _ { 0 }
^ {2 \pi} \frac {|- ( \alpha+\beta + 1 ) \overline { z } ( 1 - re ^ {i (\theta-t) } ) + ( 1 + \alpha) ( 1
- r ^ { 2 } ) e ^ { - i t } | } { |1 - re ^ {i (\theta-t) } | ^ { \alpha + \beta + 3 } }|f ( e ^ { i t }
)|^{p} \frac { d t } { 2 \pi }\right)\frac { d\theta} { 2 \pi }\\
&\leq (|c_{\alpha,\beta}|( 1 - r ^ { 2 } ) ^ {\alpha+\beta})^{p}I_6^{p}\|f\|^p_{ L ^ { p } ( \mathbb{T} )}.
\end{aligned}
\ee
Let
$$
\widetilde{C}_{\alpha,\beta}(r)=|c_{\alpha,\beta}|\int_{0}^{2\pi}\frac{|(\alpha+1)-\beta re^{is}|}
{|1+re ^ {-is}|^{-(\alpha+\beta)}}\frac { d s } { 2 \pi }.
$$
Then \eqref{eq:41018} reduces to
$$
M_{p}(r,u_{z})\leq\frac{\widetilde{C}_{\alpha,\beta}(r)}{1-r^{2}}\|f\|_{ L ^ { p } ( \mathbb{T} )}.
$$
Let
\be\label{eq:41020}
C_{\alpha,\beta}(r)=\frac{|c_{\alpha,\beta}| (|\alpha+1|+|\beta|r)}{2\pi}\int _ { 0 } ^ {2
\pi}\frac{ds}{|1+re ^ {-is}|^{-(\alpha+\beta)}}.
\ee
Then,
\be\label{eq:41021}
C_{\alpha,\beta}=|c_{\alpha,\beta}|(|\alpha+1|+|\beta|)
\frac{\Gamma(\alpha+\beta+1)}{\Gamma^{2}(\frac{\alpha+\beta+2}{2})}.
\ee
It follows that $\tilde {C} _ { \alpha,\beta } ( r ) \leq C_ { \alpha,\beta } ( r ) \leq
C_{\alpha,\beta }$. Therefore, the desired inequalities hold.

If $\alpha=\beta=0$, then $C_{0,0}( r ) = C_ {0,0}=1$.  It can be shown that the constant $C_{0,0}$ is sharp,
see \cite[Theorem 1.3]{Yong} for a  detailed proof.

The estimate for $M_{p}(r, u _ { \overline{z} } )$ is similar. We omit the proof.
\end{proof}

\section*{Declarations}

\subsection*{Ethical Approval} Not applicable.

\subsection*{Conflict of Interests}
The authors of this manuscript have no relevant financial or non-financial interests to
disclose.


\subsection*{Data Availability Statement}
Data sharing not applicable to this article as no datasets were generated or analysed during
the current study.


\subsection*{Author Contributions} All the authors of this manuscript have participated in conducting the research and writing of the manuscript.

\subsection*{Funding}
J. Qiao was supported by NSF of Hebei Science Foundation under the number A2024201013. A. Rasila was partly supported by NSF of Guangdong Province under the number 2024A1515010467 and Li Ka Shing Foundation under the number 2024LKSFG06.

\end{document}